\documentclass{birkjour}
\usepackage{amsmath, amsthm, amsfonts, mathdots}
\usepackage{amssymb}
\usepackage{amsbsy}
\usepackage{mathrsfs}

\input amssym.def
\input amssym
\usepackage{enumerate}
\usepackage{hyperref}
\hypersetup{
    pdftitle={Memory Estimation of Inverse Operators},    
    pdfauthor={Anatoly Baskakov}{Ilya Krishtal},     
    pdfsubject={Memory Estimation of Inverse Operators},   
    pdfcreator={Ilya Krishtal},   
    pdfproducer={Ilya Krishtal}, 
    pdfkeywords={Banach modules}{Beurling spectrum}{Wiener's lemma}, 
     colorlinks,%
    citecolor=blue,%
    filecolor=black,%
    linkcolor=magenta,%
    urlcolor=black
}

\chardef\bslash=`\\ 

\makeatletter
\def\verbatim{\interlinepenalty\@M \@verbatim
   \leftskip\@totalleftmargin\advance\leftskip2pc
   \frenchspacing\@vobeyspaces \@xverbatim}
\makeatother
\newtheorem{thm}{Theorem}[section]
\newtheorem{cor}[thm]{Corollary}
\newtheorem{lem}[thm]{Lemma}
\newtheorem{prop}[thm]{Proposition}

\newtheorem{ass}[thm]{Assumption}
\newtheorem{hyp}{Hypothesis}[section]

\theoremstyle{definition}
\newtheorem{defn}{Definition}[section]

\theoremstyle{remark}
\newtheorem{rem}{Remark}[section]
\newtheorem{exmp}{Example}[section]

\numberwithin{equation}{section}

\newcommand{\begeq}{\begin {equation}}
\newcommand{\eq}{\end{equation}}
\newcommand{\bs}{\begin{split}}
\newcommand{\es}{\end{split}}
\newcommand{\beh}{\begin{hyp}}
\newcommand{\eh}{\end{hyp}}
\newcommand{\bp}{\begin {prop}}
\newcommand{\ep}{\end {prop}}
\newcommand{\bt}{\begin {thm}}
\newcommand{\et}{\end {thm}}
\newcommand{\bc}{\begin {cor}}
\newcommand{\ec}{\end {cor}}
\newcommand{\bl}{\begin {lem}}
\newcommand{\el}{\end {lem}}
\newcommand{\bpf}{\begin {proof}}
\newcommand{\epf}{\end {proof}}
\newcommand{\bi}{\begin {itemize}}
\newcommand{\ei}{\end {itemize}}
\newcommand{\ben}{\begin {enumerate}}
\newcommand{\een}{\end {enumerate}}
\newcommand{\brem}{\begin {rem}}
\newcommand{\erem}{\end {rem}}

\newcommand{\bd}{\begin {defn}}
\newcommand{\ed}{\end {defn}}

\newcommand{\bex}{\begin {exmp}}
\newcommand{\eex}{\end {exmp}}

\newcommand{\la}{\langle}
\newcommand{\ra}{\rangle}

\newcommand{\F}{\mathcal{F}}

\newcommand{\TTT}{{\T\kern-.44em \T}}
\newcommand{\tTTT}{\widetilde{\T\kern-.44em \T}}

\newcommand{\ZZ}{{\mathbb Z}}

\newcommand{\RR}{{\mathbb R}}
\newcommand{\CC}{{\mathbb C}}
\newcommand{\NN}{{\mathbb N}}

\newcommand{\HH}{{\mathcal H}}

\renewcommand{\d}{\delta}
\newcommand{\s}{\sigma}
\renewcommand{\l}{\lambda}
\renewcommand{\L}{\Lambda}
\renewcommand{\a}{\alpha}

\newcommand{\M}{\mathcal{M}}

\newcommand{\T}{\mathcal{T}}

\newcommand{\JJJ}{\mathcal J}

\begin{document}

%
%
%
%
%
%
%
%
%


\title[Differential operator with an involution as a group generator]{
Linear Differential  Operator 
with an Involution as a Generator of an Operator Group} 

\author{Anatoly G. Baskakov}
\address{Department of Applied Mathematics and Mechanics \\ Voronezh State University \\ Voronezh 394693 \\ Russia}
\email{anatbaskakov@yandex.ru}
\thanks{The first author is supported in part by the Ministry of Education and Science of the Russian Federation in the frameworks of the project part of the state work quota (Project No 1.3464.2017/4.6).
The second author is supported in part by NSF grant DMS-1322127. The third author is supported in part by RFBR grant 16-01-00197.}



\author{Ilya A. Krishtal}
\address{Department of Mathematical Sciences\\ Northern Illinois University\\ DeKalb, IL 60115 \\ USA}
\email{krishtal@math.niu.edu}
\author{Natalia B. Uskova}
\address{Department of Higher Mathematics and Mathematical Physical Modeling \\ Voronezh State Technical University \\ Voronezh 394026 \\ Russia}
\email{nat-uskova@mail.ru}


\date{\today }




\begin{abstract}
We use the method of similar operators to study a mixed problem for a differential equation with an involution and an operator-valued potential function. The differential operator defined by the equation is transformed into a similar operator that is an orthogonal direct sum of simpler operators. The result is used to construct an operator group  that describes the mild solutions of the original problem. It may also serve as a justification for the use of the Fourier method to solve it.
\end{abstract}
\maketitle

\section{Introduction}

Mixed problems with an involution arise in various theoretical and applied research fields such as filtering and prediction theory \cite{KB61} and the study of subharmonic oscillations \cite{P64, P65}. Some classical geometric problems of Bernoulli and Euler may also lead to a (finite dimensional) system of differential equations with a simple involution \cite{V84}. In addition, problems with an involution are  interesting because of their relation to the problems with a Dirac operator \cite{BDS11, BK13, Mit04}. 
In this paper,
we study a mixed problem for a differential equation with an involution and an operator-valued potential function in the following form:

\begin{equation}\label{bask1}
\begin{cases}
\frac{\partial u(t, s)}{\partial t}=\frac{\partial u(t, s)}{\partial s}-\mathcal{V}(s)u(t,\omega-s),\\
u(t, 0)=u(t, \omega),\\
u(0, s)=\varphi(s),\\
t\in\mathcal{J}, \quad s\in [0, \omega],
\end{cases}
\end{equation}
where $\mathcal{J}$ is one of the intervals $(-\infty, \infty)$, $(-\infty, \beta]$, $[\alpha, \beta]$, or
$[\alpha, \infty)$. It is always assumed that $0\in\JJJ$.

We also consider a non-homogeneous problem given by
\begin{equation}\label{bask2}
\begin{cases}
\frac{\partial u(t, s)}{\partial t}=\frac{\partial u(t, s)}{\partial s}-\mathcal{V}(s)u(t,\omega-s)+f(t, s),\\
u(t, 0)=u(t, \omega),\\
u(0, s)=\varphi(s),\\
t\in\mathcal{J}, \quad s\in [0, \omega].
\end{cases}
\end{equation}

The mixed problems (\ref{bask1}) and (\ref{bask2}) and related differential operators have been studied, for example, in \cite{BKh11z, BKh11, BKh14, BKK12,  K81, SS12} in the  case  of scalar valued functions 
and with a smooth potential $\mathcal V: [0, \omega]\to \mathbb{C}$. For the homogeneous problem (\ref{bask1}), the resolvent method and contour integration were used to justify the Fourier method. The authors of \cite{BKh11z, BKh11, BKh14} also obtained results on the asymptotics of the eigenvalues and the equiconvergence of the spectral decompositions of the differential operator $L$ defined by the problem (\ref{bask1}). The spectral properties of the operator $L$ for the case of $\CC^d$-valued functions, $d\in \NN$, were  studied in \cite{BKR17}. There, the authors used the method of similar operators, which will also be the primary tool in this paper.
The method was pioneered by Friedrichs \cite{F65} and then extensively developed and used, for example,  in \cite{B86, B94, BDS11, BP17}. In this paper, we partially follow the blueprint of \cite{BKR17} to extend some of the results on the spectral properties of the operator $L$ to the case when the functions have values in an infinite dimensional Hilbert space. Our primary focus, however, is describing the group generated by the operator $L$ and studying its spectral properties.

This paper is organized as follows. In Section \ref{sect2}, we introduce the main notions and notation. In particular, in Subsection \ref{sub21}, we carefully define the mixed problems we study, their classical and mild solutions, and the differential operator $L$ that is associated with them. In Subsection \ref{sub22}, we define the notion of similarity for unbounded operators and present some known results about the spectra of similar operators. In Subsection \ref{sub23}, we state the results of this paper. The foundational result   is Theorem \ref{baskth3},  which proves that, under mild conditions on the potential function $\mathcal V$, the operator $L$ is similar to an analogous operator with a  Hilbert-Schmidt-valued potential. 
Our main contributions are Theorem \ref{baskth4}, where we describe the spectral properties of the operator $L$ and Theorem \ref{baskth6'}, where we give a more or less explicit formula for the group $\widetilde T$ that is similar to the group $T$  generated by $L$. Two other important results deal with the equiconvergence of spectral decompositions (Theorem \ref{baskth5}) and the approximation estimates for $T$ in terms of $\widetilde T$ (Theorem \ref{baskth8'}).
Section \ref{methsim} contains a brief description of the method of similar operators in an abstract setting. In Sections \ref{sect3} and \ref{sect4}, we use the method of similar operators to construct two consecutive similarity transforms for the differential operator studied in this paper. Section \ref{sect5} contains the  proofs of the results in Section \ref{sect2}. Finally, Section \ref{sect7} concludes the paper with an illustrative example that is intended to help the reader appreciate the method of similar operators, our results, and the role played by the involution in \eqref{bask1}.

\section{Main definitions and results}\label{sect2}

In this section, we provide precise definitions of the mixed problems, their various solutions, and the differential operator associated with them. We remind the reader of the relevant facts from the theory of operator semigroups and exhibit basic definitions needed for the method of similar operators. We conclude the section with the statements of the main results of this paper.

\subsection{Mixed problems and their solutions.}\label{sub21}
To make Problems \eqref{bask1} and \eqref{bask2} precise, we introduce the following notation. We let $H$ be a complex Hilbert space and $L^2 = L^2([0,\omega], H)$ be the Hilbert space of all (equivalence classes) of square summable Lebesgue measurable $H$-valued functions. The inner product on $L^2$ is given by
$$
\la x, y\ra=\frac{1}{\omega}\int_0^\omega\la x(s), y(s)\ra_H\,ds, \quad x, y\in L^2.
$$
If $H = \CC$, we shall write $L^2[0,\omega]$ instead of $L^2([0,\omega], \CC)$.

By $W_2^1=W_2^1([0, \omega], H)$ we denote the Sobolev space of continuous $L^2$ functions with derivatives
in $L^2$ and the inner product $\la x, y\ra_{W_2^1}=\la x, y\ra+\la x', y'\ra$, $x$, $y\in W_2^1$.

For  an abstract complex Hilbert space $\HH$ we denote by
$B(\HH)$  the Banach algebra of all bounded linear operators in $\HH$. 

We assume that the function $\mathcal{V}: [0, \omega]\to B(H)$ in \eqref{bask1} and \eqref{bask2} is strongly measurable \cite{DS58} and
the functions
 $s\mapsto \mathcal{V}(s)x: [0, \omega]\to H$, $x\in H$,
belong to $L^2([0, \omega], H)$. Thus, the Fourier series of the function $\mathcal{V}$ is well defined by
$$
\mathcal{V}(s)x=\sum_{n\in\mathbb{Z}}\widehat{\mathcal{V}}(n)xe^{i\frac{2\pi n}{\omega}s}, \quad s\in [0, \omega],
\quad x\in H,
$$
where the Fourier coefficients $\widehat{\mathcal{V}}(n)\in B(H)$, $n\in\mathbb{Z}$, are given by
$$
\widehat{\mathcal{V}}(n)x=\frac{1}{\omega}\int_0^\omega \mathcal{V}(s)xe^{-i\frac{2\pi n}{\omega}s}\,ds, \quad x\in H.
$$
We note that $\|\mathcal{V}(s)x\|_2^2=\sum\limits_{n\in\mathbb{Z}}\|\widehat{\mathcal{V}}(n)x\|^2<\infty$ for each $x\in H$ and $s\in [0,\omega]$. We shall, however, always assume two stronger conditions:
\begin{equation}\label{bask3}
\sum\limits_{n\in\mathbb{Z}}\|\widehat{\mathcal{V}}(n)\|^2<\infty
\end{equation}
and 
\begin{equation}\label{bask333}
\sum\limits_{j,\ell\in\mathbb{Z}}\left\|\sum_{n\ne \ell}\frac{\widehat{\mathcal{V}}(j+n)\widehat{\mathcal{V}}(\ell+n)}{\ell-n}\right\|^2<\infty.
\end{equation}
The nature of the last condition will become apparent in Section \ref{sect3}. There, in Proposition \ref{exprop}, we shall also provide sufficient conditions for \eqref{bask333} to hold.
One of those conditions implies that \eqref{bask3} and \eqref{bask333} may hold even when 
the function $s\mapsto \|\mathcal V(s)\|$ is unbounded.

By $C(\mathcal{J}, L^2)$ we shall denote the linear space of all functions $v: \mathcal{J}\times [0, \omega]\to H$ such that, for each fixed $t\in\mathcal{J}$, the function $s\mapsto v(t, s)$ belongs to $L^2=L^2([0, \omega], H)$ and the function
$$
\widetilde{v}: \mathcal{J}\to L^2, \quad (\widetilde{v}(t))(s)=v(t, s), \quad t\in\mathcal{J}, \quad s\in [0, \omega],
$$
is continuous.
If $\mathcal{J}$ is a finite interval then $C(\mathcal{J}, L^2)$ is a Banach space with the norm
 $\|v\|_\infty=\max\limits_{t\in\mathcal{J}}\|\widetilde{v}(t)\|_2$.
We call  the function $\widetilde{v}$ the \emph{associated function} to $v$ and frequently identify the two in the rest of the paper.

The function $f: \mathcal{J}\times [0, \omega]\to H$ in the non-homogeneous problem (\ref{bask2}) is assumed to belong to the space $C(\mathcal{J}, L^2)$.

Problems \eqref{bask1} and \eqref{bask2} have the following equivalent formulations in $L^2$:
\begin{equation}\label{basknew4}
\widetilde{u}_t=L\widetilde{u}, \quad \widetilde{u}(0)=\varphi,
\end{equation}
\begin{equation}\label{basknew5}
\widetilde{u}_t=L\widetilde{u}+\widetilde{f}, \quad \widetilde{u}(0)=\varphi.
\end{equation}
The operator $L: D(L)\subset L^2\to L^2$ in (\ref{basknew4}) and \eqref{basknew5}, is defined by
\begin{equation}\label{bask4}
(Ly)(s)=\frac{dy}{ds}(s)-\mathcal{V}(s)y(\omega-s), \quad s\in [0, \omega].
\end{equation}
The domain $D(L)$ is given by the periodic boundary conditions
\begin{equation}\label{bask5}
D(L)=\{y\in W_2^1: y(0)=y(\omega)\}.
\end{equation}

In the following two definitions we identify precisely two types of solutions of problems (\ref{bask1}) and (\ref{bask2}).

\bd\label{basknewdef1}(\cite{ABHN11})
By the \emph{classical solution} of Problem (\ref{bask2}) we mean a continuously differentiable function $u: \mathcal{J}\times [0, \omega]\to H$,
which belongs to the space $C(\mathcal{J}, L^2)$ and such that the associated function $\widetilde{u}: \mathcal{J}\to
L^2([0, \omega], H)$ is continuously differentiable, satisfies   $\widetilde{u}(t)\in D(L)$ for each $t\in\mathcal{J}$, and \eqref{basknew5} holds.
\ed

\bd(\cite[\S~3.1]{ABHN11})\label{basknewdef2}
A function $u\in C(\mathcal{J}, L^2)$ is a \emph{mild solution} of (\ref{basknew5}) if
$\int_0^t\widetilde{u}(s)\,ds\in D(L)$ and
$$
\widetilde{u}(t)=\varphi+L\int_0^t\widetilde{u}(s)\,ds+\int_0^t\widetilde{f}(s)\,ds, \quad t\in\mathcal{J},
$$
where the integrals are Riemann integrals of continuous functions from $\JJJ$ to  $L^2=L^2([0, \omega]; H)$.
\ed

\brem\label{baskdef1}\label{baskdef2}
We remark that the classical and mild solutions of Problem \eqref{bask1} are defined the same way if one takes $f$ and $\widetilde f$ to be the zero functions. Observe that 
a function $u: \mathcal{J}\times [0, \omega]\to H$ is a mild solution of (\ref{bask1}) if and only if 
there exists a sequence of functions $\varphi_n\in W_2^1$, $n\ge 1$, such that $\lim\limits_{n\to\infty}\varphi_n=\varphi$ in $L^2$ and
$u$ is a uniform limit on compact subsets of $\mathcal{J}\times[0,\omega]$ of a sequence of classical solutions $(u_n)$, $n\ge 1$, of Problem \eqref{bask1} with
$u_n(0,s)=\varphi_n(s)$, $s\in[0,\omega]$. 
\erem

An important tool for justifying the Fourier method for Problems \eqref{bask1} and \eqref{bask2} is provided by the following, nearly obvious, result.

\bt\label{baskth1}
Assuming \eqref{bask3} and \eqref{bask333}, the differential operator $L$ is an infinitesimal generator of a strongly continuous operator group
$
T: \mathbb{R}\to B(L^2).
$
Every classical solution $u\in C(\mathcal{J}, L^2)$ of (\ref{bask1}) is given by
\begin{equation}\label{bask6}
u(t, s)= (T(t)\varphi)(s), \quad s\in [0, \omega], \quad t\in\mathcal{J},
\end{equation}
where $\varphi\in W_2^1$ and $\varphi(0)=\varphi(\omega)$. Every mild solution is also given by (\ref{bask6}), but with $\varphi\in L^2$.
\et

The above theorem, at least in the case when the  function $s\mapsto \|\mathcal V(s)\|$ is bounded, 
follows from the general results on perturbation of operator semigroups
(see~\cite{HP57}, \cite{EN00}). We shall use the method of similar operators to present a much stronger version of Theorem~\ref{baskth1} in Theorems~\ref{baskth6'} and \ref{baskth8'}.

\brem\label{bask_rem1}
From \cite[Proposition~3.1.16]{ABHN11}, it follows that any mild solution $u\in C(\mathcal{J}, L^2)$
of (\ref{bask2}) satisfies
\begin{equation}\label{bask7}
\widetilde{u}(t)=T(t-t_0)\widetilde{u}(t_0)-\int_{t_0}^tT(t-\tau)\widetilde{f}(\tau)\,d\tau, \quad t_0, t\in\mathcal{J},
\end{equation}
where $T: \mathbb{R}\to\mathrm{End}\,L^2$ is the group of operators from Theorem~\ref{baskth1} generated by the operator $L$. Similarly, any classical solution of \eqref{bask2} satisfies \eqref{bask7} with
$\widetilde{u}(t)\in D(L)\subset W_2^1$.
\erem

The following theorem immediately follows.

\bt\label{baskth2}
Problem (\ref{bask2}) has a unique mild solution $u\in C(\mathcal{J}, L_2)$ such that
$$
\widetilde{u}(t)=T(t)\varphi+\int_0^t T(t-\tau)\widetilde{f}(\tau)\,d\tau, \quad t\in\mathcal{J}.
$$
\et

\subsection{Similar operators and direct sums.}\label{sub22}
Recall that by $\HH$ we denote an abstract complex Hilbert space. We begin with the following definition.
\bd\label{baskdef5}
Two linear operators $A_i: D(A_i)\subset\mathcal{H}\to\mathcal{H}$, $i=1, 2$, are called 
\emph{similar}, if there exists a continuously invertible operator
 $U\in B(\mathcal{H})$ such that
$$
A_1Ux=UA_2x, \quad x\in D(A_2), \quad UD(A_2)=D(A_1).
$$
The operator $U$ is called the \emph{similarity transform} of $A_1$ into $A_2$.
\ed

Directly from the above definition we have the following result about the spectral properties of similar operators.

\bl\label{basklh1}
Let $A_i: D(A_i) \subset \HH\to\HH$, $i = 1,2$, be two similar operators with the similarity transform $U$. Then the following properties hold.
\ben
\item[(1)]	We have $\s(A_1) = \s(A_2)$, $\s_p(A_1) = \s_p(A_2)$, and $\s_c(A_1) = \s_c(A_2)$, where $\s_p$ denotes the point spectrum and $\s_c$ denotes the continuous spectrum;
\item[(2)]	Assume that the operator $A_2$ admits a decomposition $A_2 = A_{21}\oplus A_{22}$ with respect to a direct sum $\HH = \HH_1\oplus \HH_2$, where $A_{21} = A_2|_{\HH_1}$ and $A_{22} = A_2|_{\HH_2}$ are the restrictions of $A_2$ to the respective subspaces. Then the operator $A_1$ admits a decomposition $A_{1} = A_{11} \oplus A_{12}$ with respect to the direct sum $\HH = \widetilde{\HH}_1\oplus \widetilde{\HH}_2$, where $A_{11} = A_1|_{\widetilde{\HH}_1}$  and   
$A_{12} = A_1|_{\widetilde{\HH}_2}$ are the restrictions of $A_1$ to the respective invariant subspaces. Moreover, if $P$ is the projection onto $\HH_1$ parallel to $\HH_2$, then $\widetilde P = UPU^{-1}$ is the  projection onto $\widetilde{\HH}_1$ parallel to $\widetilde{\HH}_2$.
\item[(3)] If $\l_0$ is an eigenvalue of the operator $A_2$ and $x$ is a corresponding eigenvector, then $y = Ux$ is  an eigenvector of the operator $A_1$ corresponding to the same eigenvalue $\l_0$.
\item[(4)] If $A_2$ is a generator of a $C_0$-semigroup (or group) $T_2: \mathbb{J}\to 
B(\mathcal{H})$, $\mathbb{J}\in\{\RR,\RR_+\}$, then the operator $A_1$ generates the $C_0$-semigroup (group)
$$
T_1(t)=UT_2(t)U^{-1}, \quad t\in\mathbb{J}, \quad T_1: \mathbb{J}\to B(\mathcal{H}),\  \mathbb{J}\in\{\RR,\RR_+\}.
$$
\een
\el

We shall need to extend Property (2) in the above lemma to the case of countable direct sums. To this end, we assume that the abstract Hilbert space $\HH$ can be written as 
\begin{equation}\label{bask8}
\mathcal{H}=\bigoplus_{\ell\in\mathbb{Z}}\mathcal{H}_\ell,
\end{equation}
where each $\mathcal{H}_\ell$, $\ell\in\mathbb{Z}$, is a closed nonzero subspace of $\HH$, $\HH_j$ is orthogonal to $\HH_\ell$ for $\ell\neq j\in\ZZ$, and each $x\in\HH$ satisfies $x=\sum\limits_{\ell\in\mathbb{Z}}x_\ell$, where
$x_\ell\in\mathcal{H}_\ell$ and $\|x\|^2=\sum\limits_{\ell\in\mathbb{Z}}\|x_\ell\|^2$. In other words, we have a disjunctive resolution of the identity 
\begeq\label{bask88}
\mathcal{P} = \{{P}_\ell,  \ell\in\mathbb{Z}\}, 
\eq
that is a system of idempotents with the following properties:

\ben
\item ${P}_\ell^*={P}_\ell$, $\ell\in\mathbb{Z}$;

\item ${P}_j{P}_\ell=\d_{j\ell} {P}_\ell$,  $j$, $\ell\in\mathbb{Z}$, where $\d_{j\ell}$ is the standard Kronecker delta;

\item The series $\sum\limits_{n\in\mathbb{Z}}{P}_\ell x$ converges unconditionally to $x\in\mathcal{H}$ and
$\|x\|^2=\sum\limits_{\ell\in\mathbb{Z}}\|{P}_\ell x\|^2$;

\item Equalities ${P}_\ell x=0$, $\ell\in\mathbb{Z}$, imply $x=0\in\HH$;

\item $\mathcal{H}_\ell=\mathrm{Im}\,{P}_\ell$, $x_\ell={P}_\ell x$, $\ell\in\mathbb{Z}$.
\een

Given a resolution of the identity $\mathcal P$ as in \eqref{bask88}, it is often convenient to represent an operator $X\in B(\HH)$ in terms of its matrix. We write such matrices as $\widetilde{X}=(X_{j\ell})_{j, \ell\in\mathbb{Z}}$, where $X_{j\ell}={P}_jX{P}_\ell$,
$j$, $\ell\in\mathbb{Z}$. Under obvious conditions on the domain, the matrix also makes sense for an unbounded operator. In the case when the matrix of a linear operator is diagonal, we get the following definition.

\bd\label{baskdef6}
We say that a closed linear operator ${A}: D({A})\subset\mathcal{H}\to\mathcal{H}$ is represented as an \emph{orthogonal direct sum} of bounded operators ${A}_\ell\in B(\mathcal{H}_\ell)$, $\ell\in\mathbb{Z}$, with respect to a decomposition
(\ref{bask8}), that is
\begin{equation}\label{bask9}
{A}=\bigoplus_{\ell\in\mathbb{Z}}{A}_\ell,
\end{equation}
if the following three properties hold.
\ben
\item $D({A})=\{x\in\mathcal{H}: \sum\limits_{\ell\in\mathbb{Z}}\|A_\ell x_\ell\|^2<\infty,
x_\ell={P}_\ell x, \ell\in\mathbb{Z}\}$ and $\mathcal{H}_\ell\subset D({A})$ for all $\ell\in\mathbb{Z}$.

\item For each $\ell\in\mathbb{Z}$ the subspace $\mathcal{H}_\ell$ is an invariant subspace of the operator ${A}$ and ${A}_\ell$ is the restriction of ${A}$ to $\mathcal{H}_\ell$. The operators
 ${A}_\ell$, $\ell\in\mathbb{Z}$, are called \emph{parts} of the operator ${A}$.

\item $Ax=\sum\limits_{\ell\in\mathbb{Z}}A_\ell x_\ell$, $x\in D({A})$, where $x_\ell={P}_\ell x$, $\ell\in\mathbb{Z}$, and the series converges unconditionally in $\HH$.
\een
\ed

We remark that $\sigma({A}_\ell)\subseteq\sigma({A})$, $\ell\in\mathbb{Z}$.
Without additional assumptions, however, the spectrum  $\sigma({A})$ may strictly contain the union of the spectra $\sigma({A}_\ell)$, $\ell\in\mathbb{Z}$, and even its closure.

\bd\label{baskdef7'}
Given a continuously invertible operator $U\in B(\HH)$ and an orthogonal decomposition \eqref{bask8}, a \emph{$U$-orthogonal decomposition} of $\HH$ is the orthogonal direct sum 
\begin{equation}\label{bask10''}
\mathcal{H}=\bigoplus_{\ell\in\mathbb{Z}}U\mathcal{H}_\ell.
\end{equation}
\ed

\bd\label{baskdef8'}
Given a continuously invertible operator $U\in B(\HH)$,
we say that a closed linear operator $A: D(A)\subset\mathcal{H}\to\mathcal{H}$ is a $U$-orthogonal direct sum of bounded linear operators $\widetilde{A}_\ell$, $\ell\in\mathbb{Z}$, with respect to decomposition
 (\ref{bask10''}), if \eqref{bask9} holds with respect to the decomposition \eqref{bask8} and $\widetilde{A}_\ell=UA_\ell U^{-1}$, $\ell\in\mathbb{Z}$. In this case, we write
$$
A= 
{\bigoplus_{\ell\in\mathbb{Z}}}\widetilde{A}_\ell.
$$
\ed
We remark that $U$-orthogonal decompositions and direct sums can be viewed as orthogonal with respect to the inner product
\[
\la x, y\ra_U=\la Ux, Uy\ra, \quad x, y\in\mathcal{H}.
\]

We shall provide an example of direct sums of operators based on the operator $L$ defined in (\ref{bask4}) and (\ref{bask5}). We write $L=L_0-V$, where $L_0: D(L_0)=D(L)\subset L^2\to L^2$ is the differential operator $L_0=\frac{d}{ds}$ and
\begin{equation}\label{baskNew13}
(Vy)(s)=\mathcal{V}(s)y(\omega-s).
\end{equation}
The operator $V$ is well defined because $D(L_0)\subset D(V)$. We shall treat the operator $L$ as the perturbation of $L_0$ by $V$.

Let $\mathcal{H}=L^2 =L^2([0, \omega], H)$. The spectrum of the operator $L_0=\frac{d}{ds}$ can be written as
$$
\sigma(L_0)=\bigcup_{\ell\in\mathbb{Z}}\sigma_\ell,
$$
where $\sigma_\ell=\{\lambda_\ell\}$, $\lambda_\ell=i\frac{2\pi \ell}{\omega}$, $\ell\in\mathbb{Z}$. The corresponding spectral projections 
${P}_\ell=P(\sigma_\ell, L_0)$, $\ell\in\mathbb{Z}$,  are given by
\begin{equation}\label{bask10'}
({P}_\ell x)(s)
e^{i\frac{2\pi \ell}{\omega}s}
=\widehat{x}(\ell)e^{i\frac{2\pi \ell}{\omega}s}, \; x\in\mathcal{H}, \; s\in [0, \omega],
\end{equation}
where $\mathcal{H}_k=\mathrm{Im}\,{P}_\ell$, $\ell\in\mathbb{Z}$, and the Fourier 
transform $\mathcal F: \HH\to \ell^2(\ZZ, H)$ is
defined by
\begeq\label{FT}
(\mathcal F x)(\ell) = \widehat{x}(\ell)=\frac{1}{\omega}\int_0^\omega x(\tau)e^{-i\frac{2\pi \ell}{\omega}\tau}\,d\tau, \quad \ell\in\mathbb{Z},\ x\in\mathcal{H}.
\eq
In particular, the operator $L_0=\frac{d}{ds}$ is an orthogonal direct sum of operators
$(L_0)_\ell=L_0|\mathcal{H}_\ell=\frac{i2\pi \ell}{\omega}I_\ell$, where $I_\ell$ is the identity operator on
$\mathcal{H}_\ell=\mathrm{Im}\,{P}_\ell$. In other words, $L_0=\bigoplus\limits_{\ell\in\mathbb{Z}} \lambda_\ell I_\ell$, $\ell\in\mathbb{Z}$.

We shall also make use of coarser decompositions of $\HH$ and $L_0$. For $m\in\mathbb{Z}_+=\mathbb{N}\cup\{0\}$, we let ${P}_{(m)}=\sum\limits_{|\ell|\leq m}{P}_\ell$ and consider a new resolution of the identity 
\begeq\label{bask99}
\mathcal{P}^{(m)} = \{{P}_{(m)}\}\cup\{{P}_\ell,  |\ell|> m\}. 
\eq
Then the operator $L_0$ may also be represented as an orthogonal direct sum
$$
L_0=(L_0)_{(m)}\oplus\bigg(\bigoplus_{|\ell|>m}(L_0)_\ell\bigg)=(L_0)_{(m)}\oplus\bigg(\bigoplus_{|\ell|>m}
\l_\ell I_\ell\bigg),
\quad m\in\mathbb{Z}_+, \quad \ell\in\mathbb{Z},
$$
where $(L_0)_{(m)}$ is the restriction of  $L_0=\frac{d}{ds}$ to
$\mathcal{H}_{(m)}=\mathrm{Im}\,{P}_{(m)}$. The above representation is with respect to the orthogonal decomposition of $L^2$ given by $L^2=\mathcal{H}_{(m)}\oplus\bigg(\bigoplus\limits_{|\ell|>m}\mathcal{H}_\ell\bigg)$.
Observe that $(L_0)_{(m)}=\bigoplus\limits_{|\ell|\leq m} \l_\ell I_\ell$ with respect to the decomposition
$\mathcal{H}_{(m)}=\bigoplus\limits_{|\ell|\leq m}\mathcal{H}_\ell$.

\subsection{Schatten-type classes and main results.}\label{sub23}

For an abstract complex Hilbert space  $\HH$, we write $\mathfrak{S}_p(\mathcal{H})$, $1\le p <\infty$, for the classical Schatten classes of compact operators in $B(\HH)$.
In particular, $\mathfrak{S}_2(\HH)$ is the
two-sided ideal of Hilbert--Schmidt operators in $B(\HH)$ with the norm
\[
\|X\|_2 = (\mathrm{tr}(XX^*))^{1/2},\quad X\in \mathfrak{S}_2(\HH),
\]
where $\mathrm{tr}(XX^*)$ is the trace of the nuclear operator $XX^*\in\mathfrak{S}_1(\HH)$.
The formula $\la X, Y\ra = \mathrm{tr}(XY^*)$, $X, Y\in\mathfrak{S}_2(\HH)$ defines an inner product in
$\mathfrak{S}_2(\HH)$. 

To formulate our results, however, we need a more general version of the Hilbert--Schmidt class.

\bd\label{HScD}
We say that an operator $X\in B(\mathcal{H})$ belongs to 
the \emph{Hilbert--Schmidt class $\mathfrak{S}_2(\mathcal{H}, \mathcal{P})$
with respect to a resolution of the identity $\mathcal{P}$} given by (\ref{bask88}),
if
\begeq\label{Fn}
\sum_{j,\ell\in\mathbb{Z}}\|{P}_j X{P}_\ell\|^2<\infty.
\eq
\ed

The norm $\|X\|_{2, \mathcal{P}}=\Big(\sum\limits_{j, \ell\in\mathbb{Z}}\|{P}_jX{P}_\ell\|^2\Big)^{1/2}$ turns $\mathfrak{S}_2(\mathcal{H}, \mathcal{P})$ into a normed linear space.
Moreover, the following two lemmas immediately follow.

\bl\label{basklh3}
The space $\mathfrak{S}_2(\mathcal{H}, \mathcal{P})$ of Hilbert--Schmidt operators with respect to (\ref{bask88}) is a Banach algebra.
\el

\bl\label{basklh5}
An operator $X\in\mathfrak{S}_2(\mathcal{H}, \mathcal{P})$ belongs to $\mathfrak{S}_2(\mathcal{H})$ if and only if
$$
\sum_{j, \ell\in\mathbb{Z}}\|{P}_jX{P}_\ell\|_2^2<\infty.
$$
\el

\brem\label{basklh4}
Observe that  if there is an $N\in\NN$ such that for each $\ell\in\ZZ$ the rank of 
${P}_\ell\in \mathcal P$ is at most $N$, then $\mathfrak{S}_2(\mathcal{H}, \mathcal{P})=\mathfrak{S}_2(\mathcal{H})$. We also note that for any $m\in\NN$ we have
$\mathfrak{S}_2(\mathcal{H}, \mathcal{P}^{(m)})=\mathfrak{S}_2(\mathcal{H},\mathcal{P})$, where the family $\mathcal P^{(m)}$ is given by \eqref{bask99}.
\erem

Before we state our main results, we remind the reader of the assumptions that we have made. The following Hypothesis \ref{hypo} applies to all of the remaining statements of this section as well as the rest of the paper, with the exception of Section \ref{methsim}.

\beh\label{hypo}
The  differential operator $L = L_0 - V: D(L)\subset L^2\to L^2$ is defined by \eqref{bask4} and \eqref{bask5}.  The operator $L_0: D(L_0)=D(L)\subset L^2\to L^2$ is the differential operator $L_0=\frac{d}{ds}$. The operator $V$ is given be \eqref{baskNew13}, where the potential function $\mathcal V$ satisfies \eqref{bask3} and \eqref{bask333}.
\eh

The foundation for the main contributions of this paper is provided in the following theorem.

\bt\label{baskth3}
Let $\HH = L^2 = L^2([0, \omega], H)$ and $\mathcal P$ be the resolution of the identity consisting of the spectral projections of the operator $L_0$.
There is a number $k\in\mathbb{Z}_+$ and a continuously invertible operator $U\in B(\HH)$ given by $U=I+W$ with
$W\in\mathfrak{S}_2(\mathcal{H}, \mathcal{P})$, such that  the operator $L=L_0-V$ is similar to the operator $L_0-V_0$, where  $V_0\in \mathfrak{S}_2(\mathcal{H}, \mathcal{P})$,
$$
LU=U(L_0-V_0),
$$
and the subspaces $\mathcal{H}_{(k)}=\mathrm{Im}\,{P}_{(k)}$
and $\mathcal{H}_\ell=\mathrm{Im}\,{P}_\ell$, $|\ell|>k$, are invariant for $V_0$.
Moreover, the operator $L$ is the $U$-orthogonal direct sum
$$
L=U\left(L_0-\left(V_{0(k)}\oplus\left(\bigoplus_{|\ell|>k}V_{0, \ell}\right)\right)\right)U^{-1}
$$
with respect to the $U$-orthogonal decomposition
$\mathcal{H}=U\mathcal{H}_{(k)}\oplus\left(\bigoplus\limits_{|\ell|>k}U\mathcal{H}_\ell\right)$.
\et

In Sections \ref{sect3}, \ref{sect4}, and \ref{sect5}, we will provide an explicit form for the operators $U$ and
$V_0$ in the above theorem. 

The first of our two main contributions is the following theorem, which describes the spectral properties of the operator $L$. We use the notation of Theorem
\ref{baskth3} and Hypothesis \ref{hypo}.

\bt\label{baskth4}
The spectrum $\s(L)$ of the operator $L$ satisfies
$$
\sigma(L)=\sigma_{(k)}\cup\left(\bigcup_{|\ell|>k}\left(\left\{i\frac{2\pi \ell}{\omega}\right\}+\sigma_\ell\right)\right),
$$
where each $\sigma_\ell=\sigma(V_{0,\ell})$, $|\ell|>k$, is the spectrum of the restriction $V_{0,\ell}$ of the operator $V_0$ to the invariant subspace $\mathcal{H}_\ell$.
Moreover, we have
$$
\sum\limits_{|\ell|>k}|\sigma_\ell|^2<\infty, \quad \mbox{ where } \quad |\sigma_\ell|=\max\limits_{\lambda\in\sigma_\ell}|\lambda|,
\quad \ell\in\mathbb{Z}.
$$
\et

We remark that the differential operators without the involution studied in \cite{BD15, BKKS17} have considerably different spectral properties.
For the finite dimensional case, a stronger version of this theorem can be found in \cite{BKR17}. We also state the following theorem that follows immediately from Theorem
\ref{baskth4} and results in \cite{B78}.

\bt
If $H$ is a finite dimensional space then each bounded mild solution $\tilde u$ of \eqref{basknew4} is a Bohr almost periodic function \cite{B78}. 
\et

The following theorem holds under the same conditions as Theorems~\ref{baskth3} and \ref{baskth4} and uses the same notation.
We also let $\widetilde{{P}}_{(k)}=P(\sigma_{(k)}, L)$ 
and $\widetilde{P}_\ell=P(\{\frac{i2\pi \ell}{\omega}\}+\sigma_\ell, L)$,
$|\ell|>k$, be the spectral projections corresponding to the sets $\sigma_{(k)}$ and $\{\frac{i2\pi \ell}{\omega}\}
+\sigma_\ell$, $|\ell|>k$, respectively. The result below establishes the equiconvergence of the spectral decompositions in the topology of $\mathfrak{S}_2(\mathcal{H}, \mathcal{P})$.
\bt\label{baskth5}
We have
$$
\lim_{n\to\infty}\left\|\widetilde{{P}}_{(k)}+\sum_{|\ell|=k+1}^n\widetilde{{P}}_\ell-{P}_{(k)}-
\sum_{|\ell|=k+1}^nP_\ell\right\|_{2, \mathcal{P}}=0.
$$
\et

The following theorem is the second main contribution of this paper.  
\bt\label{baskth6'}
Consider the differential operator $L$ defined by \eqref{bask4} and \eqref{bask5} and let $\mathcal P$ be the resolution of the identity consisting of the spectral projections of the operator $L_0$. The operator $L$
generates a $C_0$-group of operators
$$
T: \mathbb{R}\to B(\HH), \quad \HH = L^2=L^2([0, \omega], H).
$$
Moreover, there is a number $k\in\mathbb{Z}_+$ and a continuously invertible operator $U\in B(\HH)$ given by $U=I+W$ with
$W\in\mathfrak{S}_2(\mathcal{H}, \mathcal{P})$, such that the group 
$\widetilde{T}: \mathbb{R}\to B(\HH)$ given by $\widetilde T(t)=U^{-1}{T}(t)U$, $t\in\mathbb{R}$, satifies
\begeq\label{bask0}
\widetilde{T}(t)=\left(\bigoplus_{\ell=-\infty}^{-k-1}e^{t(\frac{i2\pi \ell}{\omega}I_\ell+B_\ell)}\right)\oplus e^{tB_{(k)}}\oplus\left(\bigoplus_{\ell=k+1}^\infty
e^{t(\frac{i2\pi \ell}{\omega}I_\ell+B_\ell)}\right), \quad t\in\mathbb{R}.
\eq
In the above formula we have $B_{(k)}\in B(\mathcal{H}_{(k)})$ and
$B_\ell\in B(\mathcal{H}_\ell)$, $|\ell|>k$,   with $\sum\limits_{|\ell|>k}\|B_\ell\|^2<\infty$.
\et

The number $k$ in the above theorem will be defined more explicitly  in Theorem~\ref{baskth14}. 
Observe that each operator $T(t)$, $t\in\mathbb{R}$, is a $U$-orthogonal direct sum of operators with respect to the $U$-orthogonal decomposition of $\HH = L^2$ given by~(\ref{bask10''}). The following corollary immediately follows from \eqref{bask0}. In its
formulation, we use the translation group $S:\RR\to B(\HH)$ given by  
\[
(S(t)x)(s) = x(t + s) = \sum_{\ell\in\ZZ} e^{\l_\ell t} (P_\ell x)(s),\ x\in\HH,\ t\in\RR,
\]
where we rely on the fact that $\HH=L^2([0,\omega], H)$ is isometrically isomorphic to the space $L_{2, \omega}=L_{2, \omega}(\mathbb{R}, H)$ of $\omega$-periodic $H$-valued functions on $\mathbb{R}$ that are norm-square-summable over $[0,\omega]$.

\bc\label{shiftrep}
The group $\widetilde{T}$ in \eqref{bask0} can be written as 
$$\widetilde T (t)= S(t)\Phi^r(t) = \Phi^l(t)S(t),\quad t\in\RR,$$ 
where the functions $\Phi^l, \Phi^r: \CC\to B(\HH)$ are entire functions of exponential type such that  $\Phi^l|\mathcal{H}_j = \Phi^r |\mathcal{H}_j$, $|j| > k$, and
\[
\max\{\|\Phi^l(z)\|,\|\Phi^r(z)\|\} \le e^{\beta|z|}, \ \beta = \sup\left\{\frac{2\pi k}{\omega}+\|B_{(k)}\|, \|B_j\|, |j|> k\right\}.
\]
\ec

Given a mild solution $u(t, s)=(T(t)\varphi)(s)$, $t\in\mathbb{R}$, $s\in [0, \omega]$,
$\varphi\in L^2$, of Problem (\ref{bask1}), we define its \emph{generalized Fourier series} by
\begeq\label{basksem}
\begin{split}
T(t)\varphi &=U\widetilde{T}(t)U^{-1}\varphi = U\sum_{\ell=-\infty}^{-k-1}e^{t(\frac{i2\pi \ell}{\omega}I_\ell+B_\ell)}P_\ell U^{-1}\varphi \\
& + Ue^{tB_{(k)}}P_{(k)}U^{-1}\varphi+U\sum_{\ell=k+1}^\infty e^{t(\frac{i2\pi \ell}{\omega}I_\ell+B_\ell)}P_\ell U^{-1}\varphi.
\end{split}
\eq

In the following theorem we estimate the rate of convergence of the generalized Fourier series of a mild solution of \eqref{bask1}. We let $\kappa_\ell=\sup\limits_{|n|\ge \ell}\|B_n\|$,
$\beta_\ell=\|P_\ell W\|_2$, and $\gamma_\ell=\|B_\ell\|$, $|\ell|\ge k+1$.
 Then $\lim\limits_{l\to\infty}\kappa_\ell=0$,
$\lim\limits_{|\ell|\to\infty}\beta_\ell=0$, and $(\beta_\ell)$ and $(\gamma_\ell)$ are square summable
sequences.

Using (\ref{basksem}) and Parseval's identity, we obtain the following result.
\bt\label{baskth8'}
For any function $\psi\in L^2$  we have
\begin{flalign*}
&\left\|T(t)U\psi -Ue^{tB_{(k)}}P_{(k)}\psi-\sum_{k< |\ell|\leq n}Ue^{(\frac{i2\pi \ell}{\omega}I_\ell+B_\ell)t}P_\ell\psi \right\|_2
 \\
&\leq C\left(\sum_{|\ell|\geq n+1}e^{2\gamma_\ell|t|}
(|\widehat{\psi}(\ell)|^2+\beta_\ell^2\|\psi\|^2)\right)^\frac{1}{2}\\
&\leq Ce^{\kappa_{n+1}|t|}\left(\sum_{|\ell|\geq N+1}(|\widehat{\psi}(\ell)|^2+\beta_\ell^2\|\psi\|^2)\right)^\frac{1}{2},
\quad t\in\mathbb{R},
\end{flalign*}
for each $n>k$ and some $C>0$ that is independent of $n$.
\et
We conclude the list of our results with the  following theorem that follows immediately from Theorem \ref{baskth8'}.
\bt\label{baskth8''}
The spectral and growth bounds \cite{ABHN11} of the group ${T}$ generated by the differential operator $L$ coincide.
\et

We remark that the spectral theory of abstract differential operators such as $L$ and the theory of operator (semi)groups generated by them have been extensively studied. For readers interested in the subject, we mention excellent books \cite{ABHN11, CL99, EN00} and references therein. We note that general invertibility conditions for abstract parabolic operators \cite{BK17, LMS95} do not apply for the operator $L$ we study here because $i\RR$ is not a subset of the interior of $\rho(L)$. Stability questions for the group generated by $L$ are of interest, but we leave them beyond the scope of this paper and refer to \cite{B15, T01}, among other papers on the subject.

\section{The method of similar operators.}\label{methsim}

The method of similar operators has its origins in various similarity and perturbation techniques. Among them classical perturbation methods of celestial mechanics, Ljapunov's kinematic
similarity method \cite{GKK96, Lj56, N15},  Friedrichs' method of similar operators that is used in quantum mechanics \cite{F65}, and Turner's method of similar operators  \cite{T65, U04}.

The method of similar operators that we use here originally appeared in \cite{B83, B84u}. It has many different versions that apply to various classes of differential operators \cite{B15, P16, U15, U16}. In this section, we outline the version of the method that was used in \cite{BDS11, BP17}.

The method of similar operators constructs a similarity transform for an operator $A - B:  D(A) \subset \HH \to \HH$, where the spectrum of the operator $A$ is known and has certain properties, and the operator $B$ is $A$-bounded (see Definition \ref{abdd} below). The goal of the method is to obtain an operator $B_1$ such that the operator $A-B$ is similar to $A-B_1$ and  the spectral properties of $A-B_1$ are in some sense close to those of $A$. In particular, certain spectral subspaces of $A$ will remain invariant for $A-B_1$.

\bd\label{abdd}
 Let $A: D(A) \subset \HH \to \HH$ be a linear operator. A linear operator $B: D(B) \subset \HH \to \HH$ is $A$-bounded if $D(B) \supseteq D(A)$ and $\|B\|_A = \inf\{c > 0: \|Bx\| \le c(\|x\| + \|Ax\|),\ x \in D(A)\} < \infty$.
\ed
The space $\mathfrak L_A(\HH)$ of all $A$-bounded linear operators is a Banach space with respect to the norm $\|\cdot\|_A$. Moreover, given $\l_0 \in \rho(A)$, where $\rho(A) = \CC\backslash\s(A)$ is the resolvent set of $A$, we have $B \in \mathfrak L_A(\HH)$ if and
only if $B(\l_0I - A)^{-1} \in B(\HH)$ and $\|B\|_{\l_0} = \|B(\l_0I - A)^{-1}\|_{B(\HH)}$ defines an equivalent norm in $\mathfrak L_A(\HH)$ \cite{EN00}.

The method of similar operators uses the \emph{commutator transform} 
$\mathrm{ad}_A:
D(\mathrm{ad}_A)\subset B(\mathcal{H})\to B(\mathcal{H})$ defined by
$$
\mathrm{ad}_AX = AX-XA, \quad X\in D(\mathrm{ad}_A),
$$
where the domain $D(\mathrm{ad}_A)$ contains all  $X\in B(\mathcal{H})$ such that
the following two properties hold:
\ben
\item $XD(A)\subseteq D(A)$;
\item The operator $\mathrm{ad}_AX: D(A)\to\mathcal{H}$ (uniquely) extends to a bounded operator $Y\in B(\mathcal{H})$; we then let
$\mathrm{ad}_AX=Y$.
\een

The key notion of the method of similar operators is that of an admissible triplet. Once such a triplet is constructed achieving the goal of the method becomes a  routine task.

\bd[\cite{BDS11, BP17}]\label{baskdef8}
Let $\mathcal{M}$ be a linear subspace of $\mathfrak{L}_A(\mathcal{H})$,
$J: \mathcal{M}\to\mathcal{M}$, and $\Gamma: \mathcal{M}\to B(\mathcal{H})$.
The collection $(\mathcal{M}, J, \Gamma)$ is called an \emph{admissible triplet} for the operator $A$, and the space
$\mathcal{M}$ is the \emph{space of admissible perturbations}, if the following six properties hold.

\ben
\item $\mathcal{M}$ is a Banach space that is continuously embedded in $\mathfrak{L}_A(\mathcal{H})$, i.e., $\mathcal{M}$ has a norm $\|\cdot\|_\ast$
such that there is a constant $C>0$ that yields $\|X\|_A\le  C\|X\|_\ast$ for any
$X\in\mathcal{M}$.

\item $J$ and $\Gamma$ are bounded linear operators; moreover, $J$ is an idempotent.

\item $(\Gamma X)D(A)\subset D(A)$ and
$$
(\mathrm{ad}_A\,\Gamma X)x = (X-JX)x, \quad x\in D(A), \quad   X\in\mathcal{M};
$$
moreover $Y = \Gamma X\in B(\mathcal{H})$ is the unique solution of the equation
\begin{equation}\label{bask11'}
\mathrm{ad}_A\,Y = AY-YA = X-JX,
\end{equation}
that satisfies $JY=0$.

\item $X\Gamma Y$, $(\Gamma X)Y\in\mathcal{M}$ for all $X, Y\in\mathcal{M}$, and there is a constant $\gamma>0$ such that
$$
\|\Gamma\|\le\gamma, \quad \max\{\|X\Gamma Y\|_\ast, \|(\Gamma X)Y\|_\ast\}\le \gamma\|X\|_\ast\|Y\|_\ast.
$$

\item $J((\Gamma X)JY)=0$ for all $X, Y\in\mathcal{M}$.

\item For every $X\in\mathcal{M}$ and $\varepsilon>0$ there exists a number  $\lambda_\varepsilon\in\rho(A)$,
such that $\|X(A-\lambda_\varepsilon I)^{-1}\|<\varepsilon$.
\een
\ed

To illustrate  the above definition, one should think of the operators involved in terms of infinite matrices. The operator $A$ is then represented by an infinite diagonal matrix and the operator $B$ -- by a matrix with some kind of off-diagonal decay. The transform $J$ should be thought of as a projection that picks the main (block) diagonal of an infinite matrix, whereas the transform $\Gamma$ annihilates the main (block) diagonal and weighs the remaining diagonals in accordance with equation \eqref{bask11'} thereby introducing or enhancing the off-diagonal decay. The picture will be  made more precise in Example \ref{ex2}.

To formulate the main theorem of the method of similar operators for an operator $A-B$, we use the function
$\Phi:\M\to\M$ given by
\begeq\label{bask13}
\Phi(X) = B\Gamma X-(\Gamma X)(JB)-(\Gamma X)J(B\Gamma X)+B.
\eq 
\bt[\cite{BDS11, BP17}]\label{baskth6}
Assume that $(\mathcal{M}, J, \Gamma)$ is an admissible triplet for an operator $A: D(A)\subset\mathcal{H}\to\mathcal{H}$ 
and $B\in\mathcal{M}$. Assume also that
\begeq\label{bask12}
4\gamma\|J\|\|B\|_\ast<1,
\eq
where $\gamma$ comes from the Property 4 of Definition \ref{baskdef8}. Then the 
operator $A-B$ is similar to the operator $A-JX_*$, where $X_*\in\mathcal{M}$ is the (unique) 
fixed point of the function $\Phi$ given by \eqref{bask13}, and the similarity transform of 
$A-B$ into $A-JX_*$ is given by $I+\Gamma X_*\in
B(\mathcal{H})$. Moreover, the map $\Phi:\mathcal{M}\to\mathcal{M}$ is a contraction in the ball $\{X\in\mathcal{M}: \|X-B\|_\ast\le  3\|B\|_\ast\}$, and the fixed point 
$X_*$ can be found as a limit of simple iterations: $X_0=0$, $X_1=\Phi(X_0) = B$, etc.
\et

The space $\mathcal M$ in the above theorem is typically constructed based on the  properties of the operator $B$. Condition \eqref{bask12} is there to guarantee existence
of the solution of the functional equation $\Phi(X) = X$ or, in other words, existence and uniqueness of the fixed point of $\Phi$ in the space $\mathcal M$. 

We will need the following consequence of Lemma \ref{basklh1} and Theorem \ref{baskth6}.
\bt[\cite{BP17}]\label{baskth7}
Assume that $(\mathcal{M}, J, \Gamma)$ is an admissible triplet for $A: D(A)\subset\mathcal{H}\to\mathcal{H}$, 
$B\in\mathcal{M}$ satisfies (\ref{bask12}), and $A-JX_*$ is a generator of a $C_0$-group $\widetilde{T}: \mathbb{R}\to B(\mathcal{H})$. Then the operator $A-B$ is a generator of the $C_0$-group $T: \mathbb{R}\to\mathrm{End}\,\mathcal{H}$ given by
$$
T(t)=(I+\Gamma X_*)\widetilde{T}(t)(I+\Gamma X_*)^{-1}, \quad t\in\RR,
$$
where $X_*$ is the fixed point of the function $\Phi$ in (\ref{bask13}).
\et

In this paper, we are especially interested in the case when the operator $A: D(A)\subset\mathcal{H}\to\mathcal{H}$ is the differential operator $L_0 = L+V$, where $L$ is given by \eqref{bask4} and \eqref{bask5}, and the perturbation $V$, that plays the role of the operator $B$, is given by \eqref{baskNew13}. The following spectral assumptions on the operator $A$ that are often made in the method of similar operators are clearly satisfied for the operator $L_0$.

We assume that $A$ is skew-adjoint and its spectrum
$\sigma(A)$ satisfies
$$
\sigma(A)=\bigcup_{\ell\in\mathbb{Z}}\Delta_\ell,
$$
where $\Delta_\ell$, $\ell\in\mathbb{Z}$, are compact mutually disjoint sets. By $P_\ell$, $\ell\in\mathbb{Z}$, we denote the spectral projections that correspond to 
 $\Delta_\ell$, $\ell\in\mathbb{Z}$. Then
$\mathcal P = \{P_\ell, \ell\in\mathbb{Z}\}$ is a resolution of the identity and the space $\HH$ has an orthogonal decomposition \eqref{bask8} involving the spectral subspaces $\mathcal{H}_\ell=\mathrm{Im}\,P_\ell$, $\ell\in\mathbb{Z}$. 

In this setting, it is often possible to consider an admissible triplet $(\mathcal{M}, J, \Gamma)$ such that the transform $J:\mathcal{M}\to\mathcal{M}$ is defined by
$$
JX=\sum_{\ell\in\mathbb{Z}}P_\ell XP_\ell, \quad X\in\mathcal{M}.
$$
In particular, each operator $JX$, $X\in\mathcal{M}$, is an orthogonal direct sum of operators $X_\ell=P_\ell X|\mathcal{H}_\ell$, $X_\ell\in B(\mathcal{H}_\ell)$, $\ell\in\mathbb{Z}$. We also point out that in this case the operator $A-JX$ also is an orthogonal direct sum:
\begeq\label{basknew20'}
A-JX=\bigoplus_{\ell\in\mathbb{Z}}(A_\ell-X_\ell),
\end{equation}
where $A_\ell=A|\mathcal{H}_\ell$ is the restriction of $A$ to $\mathcal{H}_\ell$.

\bt\label{baskth12}
Under the assumptions of Theorem \ref{baskth7}, the operator $A-B$ with $B\in\mathcal{M}$ is similar to the operator $A-JX_*$, $X_*\in\mathcal{M}$,
which is an orthogonal direct sum of the form (\ref{basknew20'}) with $X=X_*$ and with respect to the orthogonal decomposition  (\ref{bask8}) of $\mathcal{H}$, where $\mathcal{H}_\ell=\mathrm{Im}\,P_\ell$, $\ell\in\mathbb{Z}$.
The operator $A-B$ is then a $U$-orthogonal direct sum, where $U$ is the similarity transform of $A-B$ into $A-JX_*$.
\et

We shall denote the operator $A-JX_*$ from the above theorem by $A_0$ and its parts in decomposition \eqref{basknew20'} by $A_{0,\ell}$.

\bl\label{baskuskova_lh6}
An operator $A_0$ in (\ref{basknew20'}) is a generator of a $C_0$-group of operators
$T_0: \mathbb{R}\to B(\mathcal{H})$ if and only if
\begin{equation}\label{bask_16'}
\sup_{|t|\leq b}\sup_{\ell\in\mathbb{Z}}\|e^{tA_{0, \ell}}\|_{B(\mathcal{H}_\ell)}=C(b)<\infty,
\end{equation}
 $b\ge  1$. If (\ref{bask_16'}) holds, then the operators $T_0(t)$, $t\in\mathbb{R}$, are
 orthogonal direct sums
$$
T_0(t)=\bigoplus_{\ell\in\mathbb{Z}}e^{tA_{0, \ell}}, \quad t\in\mathbb{R},
$$
with respect to the orthogonal decomposition  (\ref{bask8}) of $\mathcal{H}$.
\el

\bpf
If (\ref{bask_16'}) holds, then the formula
$$
T_0(t)x=\sum_{\ell\in\mathbb{Z}}e^{tA_{0, \ell}}P_\ell x, \quad t\in\mathbb{R},
$$
defines an operator in $B(\HH)$. This follows from
$$
\|T_0(t)x\|^2=\sum_{\ell\in\mathbb{Z}}\|e^{tA_{0, \ell}}P_\ell x\|^2\le  C^2(b)\sum_{\ell\in\mathbb{Z}}\|P_\ell x\|^2=
C^2(b)\|x\|^2, \quad x\in\mathcal{H}.
$$
A direct computation shows that the operators $T_0(t)\in\mathrm{End}\,\mathcal{H}$, $t\in\mathbb{R}$, form a group. The group is clearly strongly continuous on the dense subset of vectors of the form
$x=\sum\limits_{|\ell|\leq n}P_\ell x$, $n\in\mathbb{Z}_+$. Therefore, it is a $C_0$-group.

The converse statement is obvious and the lemma is proved.
\epf

Thus, to study the group generated by the operator $A-B$ it suffices to  study  a group which is an orthogonal direct sum of bounded operators.

\section{The first similarity transform}\label{sect3}
In many cases, it can be difficult to define the space $\mathcal{M}$ of admissible perturbations for a given operator $A-B$. It may, however, be possible to pick a good space $\mathcal{M}$ first, and then find an operator $A-C$ that is similar to $A-B$ and such that $C\in \mathcal{M}$. In fact, this is always possible if the following assumption holds.

\begin{ass}[\cite{BP17}]\label{baskpred1}
Assume that $(\mathcal{M}, J, \Gamma)$ is  an {admissible triplet} for an operator $A$
such that the transforms $J$ and $\Gamma$ are restrictions of linear operators from 
 $\mathfrak{L}_A(\mathcal{H})$ to  $\mathfrak{L}_A(\mathcal{H})$ denoted by the same symbols. Assume also that the operator $B\in \mathfrak{L}_A(\mathcal{H})$ has the following five properties.
 
\ben
\item $\Gamma B\in B(\mathcal{H})$ and $\|\Gamma B\|<1$;

\item $(\Gamma B)D(A)\subset D(A)$;

\item $B\Gamma B$, $(\Gamma B)JB\in\mathcal{M}$;

\item $A(\Gamma B)x-(\Gamma B)Ax=Bx-(JB)x$, $x\in D(A)$;

\item For any  $\varepsilon>0$ there is $\lambda_\varepsilon\in\rho(A)$ such that
$\|B(A-\lambda_{\varepsilon}I)^{-1}\|<\varepsilon$.
\een
\end{ass}
\bt[\cite{BP17}]\label{baskth8}
If Assumption \ref{baskpred1} holds then the operator $A-B$ is similar to   $A-JB-B_0$, where
$B_0=(I+\Gamma B)^{-1}(B\Gamma B-(\Gamma B)JB)$. 
The similarity transform is given by $I+\Gamma B$ so that
$$
(A-B)(I+\Gamma B)=(I+\Gamma B)(A-JB-B_0).
$$
\et

In this section, we use the above theorem to find an operator $L_0-\widetilde V$ that is similar to the operator $L = L_0 - V$ given by \eqref{bask4} and \eqref{bask5}, see also \eqref{baskNew13}. Hereinafter, we let $\mathcal{H}=L^2=L^2([0, \omega]; H)$. We also use the fact that $\HH$ is isometrically isomorphic to the space $L_{2, \omega}=L_{2, \omega}(\mathbb{R}, H)$ of $\omega$-periodic $H$-valued functions on $\mathbb{R}$ that are norm-square-summable over $[0,\omega]$.

We  let $\mathcal{P}$ be the resolution of the identity consisting of the Riesz projections
$P_\ell = P(\{i2\pi\ell/\omega\}, L_0)$, $\ell\in\mathbb{Z}$, given by (\ref{bask10'}), and choose 
$\mathcal M = \mathfrak{S}_2(\mathcal{H}, \mathcal{P})$ -- the Hilbert--Schmidt class with respect to $\mathcal P$ (see Definition \ref{HScD}).

Recall that an operator $X: D(L_0)\subset\mathcal{H}\to\mathcal{H}$, $X\in\mathfrak{L}_{L_0}(\mathcal{H})$, can be represented via its matrix
 $(X_{j\ell})$, $j, \ell\in\mathbb{Z}$, where $X_{j\ell}=
{P}_jX{P}_\ell\in B(\mathcal{H})$, $j, \ell\in\mathbb{Z}$.

\bl\label{basknewlh1}
Assume that  $X\in\mathfrak{L}_{L_0}(\mathcal{H})$ satisfies \eqref{Fn}. Then $X$ can be uniquely extended to a bounded operator in
$B(\mathcal{H})$, denoted by the same symbol, so that $X\in\mathfrak{S}_2(\mathcal{H}, \mathcal{P})$ and $\|X\|\leq \|X\|_{2, \mathcal{P}}$.
\el

\bpf
Consider $D_f = \{x \in \HH: P_\ell x = 0$ for all but finitely many $\ell\in\ZZ\}$. Clearly, 
$D_f\subseteq D(L_0)\subseteq D(X)$ and $D_f$ is dense in $\HH$. Given $x\in D_f$, we have
\[
\|Xx\|^2 = \sum_{j\in\mathbb{Z}}\|P_jXx\|^2 = \sum_{j\in\mathbb{Z}}
\left\|\sum_{\ell\in\mathbb{Z}}P_jXP_\ell x\right\|^2\le \sum_{j\in\mathbb{Z}}
\left(\sum_{\ell\in\mathbb{Z}}\|P_jXP_\ell\|\|P_\ell x\|\right)^2.
\]
Using Schwartz's inequality, we get 
\[
\|Xx\|^2 \le \sum_{j\in\mathbb{Z}}
\sum_{\ell\in\mathbb{Z}}\|P_jXP_\ell\|^2\sum_{\ell\in\ZZ}\|P_\ell x\|^2
=\sum_{j\in\mathbb{Z}}
\sum_{\ell\in\mathbb{Z}}\|P_jXP_\ell\|^2\|x\|^2.
\]
Since $D_f$ is dense in $\HH$, the operator $X$ extends to an operator in $B(\HH)$ and the extension satisfies $X\in\mathfrak{S}_2(\mathcal{H}, \mathcal{P})$ with the required estimate.
\epf

Thus, we can determine if an operator in $\mathfrak{L}_{L_0}(\mathcal{H})$ belongs to
$\mathfrak{S}_2(\mathcal{H})$ solely from the estimate \eqref{Fn} for its matrix.

Presently, we proceed to define the transforms $J, \Gamma: \mathfrak{L}_{L_0}(\mathcal{H}) \to
\mathfrak{L}_{L_0}(\mathcal{H})$ that can be used in Theorem \ref{baskth8}. These transforms will have two equivalent representations: an integral one and a matrix one. For the integral representation we use the translation representation $S:\RR\to \HH$ mentioned before:
\[
(S(t)x)(s) = x(t + s) = \sum_{\ell\in\ZZ} e^{\l_\ell t} (P_\ell x)(s),\ x\in\HH,\ t\in\RR.
\]
For an operator $X \in B(\HH)$  with the matrix entries $X_{j\ell} = P_jXP_\ell$, $j,\ell\in
\ZZ$,  we have
\begeq\label{xk}
X_nx = \sum_{j-\ell=n} X_{j\ell}x = \frac1\omega\int_0^\omega e^{-\l_n t} S(t)XS(-t)xdt,\ x\in\HH,
\eq
for the $n$-th diagonal of $X$.

\bex
Computing the matrix form of the operator $V$ in \eqref{baskNew13} we get
\[(V_{j\ell}x)(t)=(P_jVP_\ell x)(t) = (VP_\ell x)\widehat{\ }(j)e^{\l_j t} = \widehat{\mathcal{V}}(j+\ell)\hat x(\ell)e^{\l_j t}
\]
so that $V_{j\ell} = \F^*\widehat{\mathcal{V}}(j+\ell)\F$, where $\F$ is the Fourier transform defined by \eqref{FT}. In other words, the matrix of the operator $V$ can be written
in the following equivalent form:
 \begeq\label{Vmatr}
V \sim\left(
\begin{array}{ccccc}
\ddots & \ddots & \vdots & \iddots & \iddots \\
\cdots & \widehat{\mathcal{V}}(-2) &\widehat{\mathcal{V}}(-1)  & \widehat{\mathcal{V}}(0)  & \cdots \\
\cdots & \widehat{\mathcal{V}}(-1) & \widehat{\mathcal{V}}(0) & \widehat{\mathcal{V}}(1) & \cdots \\
\cdots & \widehat{\mathcal{V}}(0) & \widehat{\mathcal{V}}(1) & \widehat{\mathcal{V}}(2) & \cdots \\
\iddots & \iddots & \vdots & \ddots & \ddots 
\end{array}
\right).
\eq
\eex

For $X \in B(\HH)$, the transforms $J$ and $\Gamma$ are defined by
\begeq\label{Jtrans}
(JX)x = X_0x = \frac1\omega\int_0^\omega  S(t)XS(-t)xdt, \ x\in\HH,
\eq
and
\begeq\label{Gtrans}
(\Gamma X)x =  \frac1\omega\int_0^\omega f(t) S(t)XS(-t)xdt, \ x\in\HH,
\eq
where $f: \RR \to \CC$ is an $\omega$-periodic function given by
\begeq\label{fgamma}
f(t) = \frac\omega 2 -t \sim\sum\limits_{n\neq 0}\frac\omega{2\pi ni}  {e^{i\frac{2\pi n}\omega t}},\ t\in(0,\omega).
\eq
We deduce from \eqref{xk} that 
\begeq\label{Gcoef}
(\Gamma X)_0 = 0\ \mbox{and} \ (\Gamma X)_n = \frac\omega{2\pi ni}X_n= \frac1{\l_n}X_n.
\eq

Similarly, for each $m \in \ZZ_+$ and $X\in B(\HH)$, we define transforms $J_m$ and $\Gamma_m$ via
\begeq\label{Jm}
J_mX = P_{(m)}XP_{(m)}+\sum_{|\ell|>m} P_\ell XP_\ell = JX-P_{(m)}(JX)P_{(m)}+P_{(m)}XP_{(m)}
\eq
and
\begeq\label{Gm}
\Gamma_mX =  \Gamma X-P_{(m)}(\Gamma X)P_{(m)} =\sum_{\substack{\max\{|j|,|\ell|\}> m \\ {j\neq \ell}}}\frac{P_jXP_\ell}{\l_j-\l_\ell}, %
\eq
so that $J_0 = J$ and $\Gamma_0 = \Gamma$.


Next, we extend the definitions of the transforms so that $J_m$, $\Gamma_m$:  	
$\mathfrak L_{L_0}(\HH) \to \mathfrak L_{L_0}(\HH)$, $m \in \ZZ_+$. Given $\l_0 \in \rho({L_0})$, we define
\begeq\label{Jext}
J_mX = J_m(X({L_0}-\l_0 I)^{-1}) ({L_0}-\l_0 I), \quad X\in  \mathfrak L_{L_0}(\HH),
\eq
and
\begeq\label{Gext}
\Gamma_mX = \Gamma_m(X({L_0}-\l_0 I)^{-1}) ({L_0}-\l_0 I), \quad X\in  \mathfrak L_{L_0}(\HH).
\eq
We observe that these extensions do not depend on the choice of $\l_0\in \rho({L_0})$ and the formulas \eqref{Jext} and \eqref{Gext} still hold. Moreover, if $x \in D({L_0})$, the formulas \eqref{Jtrans} and \eqref{Gtrans} also remain valid.
If the operators $J_mX$ and $\Gamma_mX$ admit continuous extensions, we shall denote these extensions by the same symbols and write $J_mX$, $\Gamma_mX \in B(\HH)$.

\bex\label{ex2}
Computing $J$ and $\Gamma$ transforms for the operator $V$ in \eqref{baskNew13} we get 
\[
\begin{split}
(JVx)(s) &= \frac1\omega\int_0^\omega (S(t)VS(-t)x)(s)dt = \frac1\omega\int_0^\omega (S(t)Vx)(s-t)dt \\
& =\frac1\omega\int_0^\omega \mathcal V(s+t)x(\omega - s - 2t)dt = \frac1{2\omega}\int_0^\omega \mathcal V\left(\frac{s+\omega-\tau}2\right)x(\tau)d\tau
\end{split}
\]
and, similarly,
\[
\begin{split}
(\Gamma Vx)(s) &= \frac1\omega\int_0^\omega f(t)(S(t)VS(-t)x)(s)dt\\
&   = \frac1{2\omega}\int_0^\omega f\left(\frac{\omega-s-\tau}2\right)\mathcal V\left(\frac{s+\omega-\tau}2\right)x(\tau)d\tau,
\end{split}
\]
where $f$ is given by \eqref{fgamma}. Thus, $JV$ and $\Gamma V$ are, indeed, integral operators. Moreover,
we also get that their matrices have the following equivalent form:
\begeq\label{JVmatr}
JV \sim\left(
\begin{array}{ccccc}
\ddots & \ddots & \vdots & \iddots & \iddots \\
\cdots & \widehat{\mathcal{V}}(-2) &0  & 0  & \cdots \\
\cdots & 0 & \widehat{\mathcal{V}}(0) &0 & \cdots \\
\cdots & 0 & 0 & \widehat{\mathcal{V}}(2) & \cdots \\
\iddots & \iddots & \vdots & \ddots & \ddots 
\end{array}
\right)
\eq 
and
\begeq\label{GVmatr}
\Gamma V \sim \frac\omega{2\pi}\left(
\begin{array}{cccccc}
\ddots & \ddots & \vdots & \vdots& \iddots & \iddots \\
\cdots & 0 &-\widehat{\mathcal{V}}(-1)  & -\frac12\widehat{\mathcal{V}}(0) & -\frac13\widehat{\mathcal{V}}(1)  & \cdots \\
\cdots & \widehat{\mathcal{V}}(-1) & 0 & -\widehat{\mathcal{V}}(1)& -\frac12\widehat{\mathcal{V}}(2) &\cdots \\
\cdots & \frac12\widehat{\mathcal{V}}(0) & \widehat{\mathcal{V}}(1) & 0 & -\widehat{\mathcal{V}}(3) & \cdots \\
\cdots & \frac13\widehat{\mathcal{V}}(1) & \frac12\widehat{\mathcal{V}}(2)  & \widehat{\mathcal{V}}(3) & 0& \cdots \\
\iddots & \iddots & \vdots & \vdots& \ddots & \ddots 
\end{array}
\right).
\eq
Using \eqref{bask3} and Lemma \ref{basknewlh1}, we observe  that the operators $JV$ and $\Gamma V$ belong to the class $\mathcal M = {\mathfrak S}_2(\HH,\mathcal P)$ of Hilbert--Schmidt operators with respect to $\mathcal P$. 

To use Assumption \ref{baskpred1}, we need the operator $Z = V\Gamma V$. Computing its matrix, we get 
\begeq\label{wmn}
Z_{j\ell} = \F^* \left(\sum_{n\neq \ell}\frac1{\ell-n}\widehat{\mathcal{V}}(j+n) \widehat{\mathcal{V}}(\ell+n)\right)\F,
\eq
where $\F$ is the Fourier transform given by \eqref{FT}.
It follows from \eqref{bask333} and Lemma \ref{basknewlh1} that $Z\in {\mathfrak S}_2(\HH, \mathcal P)$. 

\eex

The above example together with \eqref{Jm} and \eqref{Gm}  leads to the following result.

\bl\label{basklh7}
The operators $J_mV$, $\Gamma_mV$ and $Z_m = V\Gamma_mV$, $m\in\ZZ_+$, all belong to $\mathfrak{S}_2(\mathcal{H}, \mathcal{P})$.
\el

Since the condition \eqref{bask333} is rather nebulous, we provide a few natural sufficient conditions for it in the following proposition. To prove one of them, 
we use the integral representation of the operator $Z= V\Gamma V$:
\begeq\label{VGV}
(Zy)(s) = \frac1{2\omega}\int_0^{2\omega}f\left(\frac{s-\tau}2\right)\mathcal V\left(\frac{2\omega-s-\tau}2\right)\mathcal V(s)y(\tau)d\tau,
\eq
where $f$ is given by \eqref{fgamma}.
\begin{prop}\label{exprop}
Assume that one (or more) of the following assumptions hold.
\ben
\item $\sum\limits_{n\in\ZZ} \|\widehat{\mathcal V}(n)\|< \infty$.
\item $\mathcal V \in L^2([0,\omega], \mathfrak S_2(H))$.
\een
Then \eqref{bask333} holds. 
\end{prop}

\bpf
Assuming the first condition, we have
\[
\sum_{j,\ell\in\ZZ}\left\|\sum_{n\neq \ell}\frac1{\ell-n}\widehat{\mathcal{V}}(j+n) \widehat{\mathcal{V}}
(\ell+n)\right\|^2 = \sum_{j,\ell\in\ZZ}\left\|\sum_{n\neq 0}\frac1{n}\widehat{\mathcal{V}}(j+n) \widehat{\mathcal{V}} (2\ell+n)\right\|^2
\]
\[
\le \sum_{j,\ell\in\ZZ}\left(\sum_{r\neq 0}\|\widehat{\mathcal{V}}(j+r)\|\right)
\left(\sum_{n\neq 0}\frac1{|n|^2}\|\widehat{\mathcal{V}}(j+n)\| \|\widehat{\mathcal{V}}
(2\ell+n)\|^2\right) < \infty,
\]
where we used Schwartz's inequality.

Assuming the second condition, it follows from \eqref{VGV} and
\[
\int_0^\omega\int_0^\omega \left\|\mathcal V\left(\omega-s-\tau\right)\mathcal V(2s) \right\|_2^2
dsd\tau< \infty
\]
that $Z\in \mathfrak S_2(\HH)$ so that \eqref{bask333} follows from \eqref{wmn}.
\epf

In the following lemma we consider a few more properties from the Assumption \ref{baskpred1}.
\bl\label{basklh8}
The operators $\Gamma_mV$, $m\in\mathbb{Z}_+$, have the following properties.
\ben
\item $(\Gamma_mV)D(L_0)\subset D(L_0)$;

\item $L_0(\Gamma_mV)x-(\Gamma_mV)L_0x = (V-J_mV)x$, $x\in D(L_0)$;

\item For each $\varepsilon>0$ there is $\lambda_\varepsilon\in \rho(L_0)$ such that
$\|V(L_0-\lambda_\varepsilon I)^{-1}\|_{2,\mathcal P}<\varepsilon$.
\een
\el
\bpf
From Lemma \ref{basknewlh1}, we have $(\Gamma V)\HH_\ell \subseteq D(L_0)$, $\ell\in\ZZ$, and
 \eqref{JVmatr} and \eqref{GVmatr} immediately imply
\[
L_0(\Gamma V)P_\ell = (\Gamma V)L_0 P_\ell + (V - JV)P_\ell, \ \ell\in\ZZ.
\]
From \eqref{Jm} and \eqref{Gm} we deduce
\[
L_0(\Gamma_m V)P_\ell = (\Gamma_m V)L_0 P_\ell + (V - J_mV)P_\ell, \ \ell\in\ZZ.
\]
It follows that for any $n\in\NN$ we have
\begin{flalign*}
L_0P_{(n)}(\Gamma_mV)(L_0-\lambda_\varepsilon I)^{-1}y &=
P_{(n)}(\Gamma_mV)L_0(L_0-\lambda_\varepsilon I)^{-1}y + \\
&+ P_{(n)}(V-J_mV)(L_0-\lambda_\varepsilon I)^{-1}y=P_{(n)}\mathcal{L}y,
\end{flalign*}
where $\mathcal{L}$ is given by
$$
\mathcal{L}y=(\Gamma_mV)L_0(L_0-\lambda_\varepsilon I)^{-1}y +
(V-J_mV)(L_0-\lambda_\varepsilon I)^{-1}y.
$$
Since $\mathcal L\in B(\HH)$, we have $P_{(n)}\mathcal L y\to \mathcal L y$ for all $y\in\HH$.
Using the fact that the operator $L_0$ is closed, we deduce Properties {\it{1}} and {\it{2}}.

To prove Property {\it{3}}, observe that the matrix $Y_{j\ell}$ elements of the operator $Y = V(L_0-\lambda_\varepsilon I)^{-1}$  satisfy $Y_{j\ell} = \F^*\frac{\widehat{\mathcal{V}}(j+\ell)}{\l_j-\l_\varepsilon}\F$. It follows that 
\[
\|Y\|_{2,\mathcal P}^2 = \sum_{j,\ell\in\ZZ} \left\|\frac{\widehat{\mathcal{V}}(j+\ell)}{\l_j-\l_\varepsilon}\right\|^2 = \frac{\omega^2}{4\pi^2} \sum_{\ell\in\ZZ} \left\|{\widehat{\mathcal{V}}(\ell)}\right\|^2 \sum_{j\in\ZZ} \frac{1}{j^2+n^2} < \varepsilon
\]
for $\lambda_\varepsilon = i\lambda_n$ with a sufficiently large $n$.
\epf

The final condition in the Assumption \ref{baskpred1} is proved in the following lemma.

\bl\label{basklh9}
There is $m \ge 0$ such that $\|\Gamma_m V\|_{2,\mathcal P} < 1$.
\el
\bpf
Since $\Gamma V \in \mathfrak S_2(\HH,\mathcal P)$, we have
$$\lim\limits_{m\to\infty}\|\Gamma_mV\|_{2,\mathcal P}=\lim\limits_{m\to\infty}\|\Gamma V-P_{(m)}(\Gamma V)P_{(m)}\|_{2,\mathcal P}=0$$ and the result follows.
\epf

Applying Theorem \ref{baskth8} and the above lemma, we obtain the main result of this section.

\bt\label{baskthsimilar}
There exists $m\in\mathbb{Z}_+$ such that the operator $I+\Gamma_mV$ is invertible and
$(I+\Gamma_mV)^{-1} - I\in \mathfrak S_2(\HH,\mathcal P)$. Moreover, the operator
$L_0-V$ is similar to the operatpr $L_0-\widetilde{V}$, where
\begin{equation}\label{basktildev}
\widetilde{V}=J_mV+(I+\Gamma_mV)^{-1}(V\Gamma_mV-(\Gamma_mV)J_mV)\in\mathfrak{S}_2(\mathcal{H}, \mathcal{P}).
\end{equation}
More precisely,
$
(L_0-V)(I+\Gamma_mV)=(I+\Gamma_mV)(L_0-\widetilde{V}).
$
\et

\section{The second similarity transform}\label{sect4}

In this section, we construct another similarity transform that will allow us to prove the main results. We begin with constructing a space $\mathcal M$ of admissible perturbations that is slightly smaller than $\mathfrak{S}_2(\mathcal{H}, \mathcal{P})$.

For any $X\in\mathfrak{S}_2(\mathcal{H},
\mathcal{P})$ we define 
\begin{equation}\label{baskalpha}
\alpha_n(X)=\|X\|_2^{-\frac{1}{2}}\max\left\{\left(\sum\limits_{\substack{|\ell|\ge n \\
\ell\in\mathbb{Z}}}\|{P}_\ell X\|_{2, \mathcal{P}}^2\right)^\frac{1}{4}, \left(\sum\limits_{\substack{|\ell|\ge n \\
\ell\in\mathbb{Z}}}\|X{P}_\ell\|_{2, \mathcal{P}}^2\right)^\frac{1}{4}\right\}, \, n\in\mathbb{Z}.
\end{equation}
The sequence $(\alpha_n(X))$, $n\in\mathbb{Z}$, has the following properties.

\ben
\item $\alpha_n(X)=\alpha_{-n}(X)$, $n\in\mathbb{Z}$.

\item $\lim\limits_{|n|\to\infty}\alpha_n(X)=0$. 

\item $\alpha_n(X)\le 1$ for each $n\in\mathbb{Z}$.

\item $\alpha_n(X)\ge\alpha_{n+1}(X)$, $n\ge 0$.

\item $\alpha_n(X)\ne 0$ for all $n\in\mathbb{Z}$, provided that $P_{(m)}XP_{(m)}\ne X$ for all $m\in\mathbb{Z}_+$.

\item We have
$$
\sum_{\substack{n\in\mathbb{Z}\\ \alpha_n(X)\neq 0}}\frac{\|X{P}_n\|_{2, \mathcal{P}}^2+\|{P}_nX\|_{2, \mathcal{P}}^2}{(\alpha_n(X))^2} <\infty.
$$
\een

For any $X\in\mathfrak{S}_2(\mathcal{H}, \mathcal{P})$ we define a self-adjoint operator $F_X$ by
$$
F_X=\sum_{n\in\mathbb{Z}}\alpha_n(X){P}_n.
$$
Observe that $F_X\in B(\mathcal{H})$ can be viewed as a function of the skew-adjoint operator $L_0$:
$F_X=f_X(L_0)$, where $f_X: \sigma(L_0)\to\mathbb{R}_+$, $f_X(\lambda_n)=\alpha_n(X)$, $n\in\mathbb{Z}$, and
$\|F_X\|=\max\limits_{n\in\mathbb{Z}}|\alpha_n(X)|=1$.

We are interested in the case when $X$ is  the perturbation $\widetilde{V}\in\mathfrak{S}_2(\mathcal{H}, \mathcal{P})$
given by (\ref{basktildev}). We will assume that
 $P_{(n)}\widetilde{V}P_{(n)}\ne \widetilde{V}$ for all $n\in\mathbb{Z}_+$. Otherwise, $L_0-\widetilde{V}$ can be represented for some $n\in\NN$ as
 an orthogonal direct sum of an operator $(L_0-\widetilde{V})|{\mathcal{H}_{(n)}}$ and the operator $L_0|{\mathcal{H}^\perp_{(n)}}$, where $\mathcal{H}^\perp_{(n)}$ is the orthogonal complement of $\mathcal{H}_{(n)}$. In this case, the remaining proofs are trivial.
 To simplify the notation, we shall write $\a_n$ and $F$ instead of $\a_n(\widetilde V)$ and $F_{\widetilde{V}}$, respectively. 

We let $\mathcal M$ be the space of all operators $X\in\mathfrak S_2(\HH,\mathcal P)$ such that
$$
X=X_l F \quad\mbox{and}\quad X=FX_r,
$$
where $X_l, X_r\in\mathfrak{S}_2(\mathcal{H}, \mathcal{P})$. We let $\|X\|_{\mathcal{M}}=\max\{\|X_l\|_{2, \mathcal{P}},
\|X_r\|_{2, \mathcal{P}}\}$. Clearly,  $\|X\|_2\le \|X\|_{\mathcal{M}}$, $X\in\mathcal{M}$.

Our assumption on $\widetilde{V}$ implies that $\a_n(\widetilde{V})\neq 0$ for all $n\in\ZZ$. Therefore, $\ker F\neq\{0\}$ and the space $\mathcal M$ is a Banach space.

We remark that for any  $X\in\mathfrak{S}_2(\mathcal{H}, \mathcal{P})$ we have
$$
X=\left(\sum_{n\in\mathbb{Z}}\frac{1}{\alpha_n(X)}X{P}_n\right)F_X=
F_X\left(\sum_{n\in\mathbb{Z}}\frac{1}{\alpha_n(X)}{P}_nX\right).
$$
It follows that  $\widetilde{V}\in \mathcal M$.

Next, we observe that $\mathcal M$ is invariant for the transforms
$J_k$ and $\Gamma_k$,
$k\ge0$, defined by \eqref{Jm} and \eqref{Gm}. Moreover, for each $k\ge 0$, we have
$$
J_k(X_l F)=(J_kX_l)F, \quad J_k(FX_r)=F(J_kX_r),
$$
$$
\Gamma_k(X_l F)=(\Gamma_kX_l)F, \quad \Gamma_k(FX_r)=F(\Gamma_kX_r),
$$
where $X_r$, $X_l\in\mathfrak{S}_2(\mathcal{H}, \mathcal{P})$.

To estimate $\|\Gamma_k(XF)\|_{2, \mathcal{P}}$ and $\|\Gamma_k(FX)\|_{2, \mathcal{P}}$, $X\in\mathfrak{S}_2(\mathcal{H},
\mathcal{P})$, we introduce two sequences $(\alpha_n')$, $n\in\mathbb{N}$ and $(\widetilde{\alpha}_n)$, $n\in\mathbb{N}$,
given by
$$
\alpha_n'=\max\limits_{\substack{|i|\ge  n \\ |j|<n}}\frac{|\alpha_i-\alpha_j|}{|i-j|}, \quad
\widetilde{\alpha}_n=\frac{\omega}{2\pi}(2\alpha_n+\alpha_n'), \quad n\in\mathbb{N}.
$$
It is clear that $(\alpha_n')$ and $(\widetilde{\alpha}_n)$ belong to the space $c_0(\NN)$
of sequences vanishing at infinity.

The following  lemma is an analog of \cite[Lemma~3]{BDS11} and \cite[Lemma~3.1]{BKR17}.
\bl\label{basklh12}
For any $k\in\mathbb{Z}_+$ and $X\in\mathfrak{S}_2(\mathcal{H}, \mathcal{P})$ we have
$$
\max\{\|\Gamma_k(XF)\|_{2, \mathcal{P}}, \|\Gamma_k(FX)\|_{2, \mathcal{P}}\}\leq \widetilde{\alpha}_{k+1}\|X\|_{2, \mathcal{P}}.
$$
\el
\bpf
Let ${P}^{(k)}=I-{P}_{(k)}$, $k\in\mathbb{Z}_+$. Then
$$\|F{P}^{(k)}\|=
\left\|\sum\limits_{n\in\mathbb{Z}}\alpha_n{P}_n{P}^{(k)}\right\|\le \alpha_{k+1}.$$ From the definition of the transforms $\Gamma_k$, it follows that
\begeq\label{G1} 
\bs
\|\Gamma_k(XF)\|_{2, \mathcal{P}}&=\|\Gamma(XF{P}^{(k)}) + \Gamma({P}^{(k)}XF{P}_{(k)})\|_{2, \mathcal{P}}
 \\
&\le \frac{\omega}{2\pi}\alpha_{k+1}\|X\|_{2, \mathcal{P}}+\|\Gamma({P}^{(k)}XF{P}_{(k)})\|_{2, \mathcal{P}}.
\end{split}
\eq
We write $\Gamma({P}^{(k)}XF{P}_{(k)})$ as
$$
\Gamma({P}^{(k)}XF{P}_{(k)})=\Gamma({P}^{(k)}FX{P}_{(k)})+
\Gamma({P}^{(k)}(XF-FX){P}_{(k)}).
$$
The operator matrix  entries $({P}_\ell Z_{(k)}\mathcal{P}_j)$, $\ell, j\in\mathbb{Z}$, 
of the operators $Z_{(k)}=\Gamma({P}^{(k)}(XF-FX){P}_{(k)})$, $k\ge 0$, satisfy
$$
{P}_\ell Z_{(k)}{P}_j=\frac{f(\lambda_\ell)-f(\lambda_j)}{\lambda_\ell-\lambda_j}{P}_\ell X{P}_j=
\frac{\omega}{2\pi i}\frac{\alpha_\ell-\alpha_j}{\ell-j}{P}_\ell X{P}_j,
$$
where $|\ell|\ge k+1$, $|j|\leq k$, and ${P}_\ell Z_{(k)}{P}_j=0$ otherwise. Therefore,
\begeq\label{G2}
\bs
\|\Gamma({P}^{(k)}XF{P}_{(k)})\|_{2, \mathcal{P}} &\le  \frac{\omega}{2\pi}\alpha_{k+1}\|X\|_{2, \mathcal{P}}
+\frac{\omega}{2\pi}\max\limits_{\substack{|\ell|\ge  k+1 \\ |j|\leq  k}}\frac{|\alpha_\ell-\alpha_j|}{|\ell-j|}
\|X\|_{2, \mathcal{P}}\\
&= \frac{\omega}{2\pi}(\alpha_{k+1}+\alpha_{k+1}')\|X\|_{2, \mathcal{P}}.
\end{split}
\eq

From \eqref{G1} and \eqref{G2} we obtain
$
\|\Gamma_k(XF)\|_{2, \mathcal{P}} \le \widetilde{\alpha}_{k+1}\|X\|_{2, \mathcal{P}}.
$
The estimate for the norm of  $\Gamma_k(FX)$, $X\in\mathfrak{S}_2(\mathcal{H}, \mathcal{P})$,
is obtained in a similar fashion.
\epf

The next  lemma is an analog of \cite[Lemma~4]{BDS11} and \cite[Lemma~3.2]{BKR17}.

\bl\label{basklh13}
The collection $(\mathcal{M}, J_k, \Gamma_k)$, $k\in\ZZ_+$, is an admissible triplet for the operator $L_0$, and the constant
 $\gamma=\gamma_k$ in Definition~\ref{baskdef8} satisfies
$
\gamma_k\leq \widetilde{\alpha}_{k+1}. 
$
\el
\bpf
As we observed before, $\mathcal{M}$ is a Banach space that is continuously embedded into
$\mathfrak{S}_2(\mathcal{H}, \mathcal{P})$. It follows that $\mathcal{M}$ is continuously embedded in $\mathfrak{L}_{L_0}(\mathcal{H})$. Therefore, Property {1}
of Definition~\ref{baskdef8} holds.

Properties 2 and 5 follow immediately from the definitions of $J_k$ and $\Gamma_k$, $k\ge 0$, and Lemma~\ref{basklh5}.

Properties 3 and 6 were proved in Lemma~\ref{basklh8}.

It remains to prove Property 4. Let $X=X_{l}F\in\mathcal{M}$ and $Y=Y_{l}F\in\mathcal{M}$, where $X_l$ end
$Y_l\in\mathfrak{S}_2(\mathcal{H}, \mathcal{P})$. Then $X\Gamma_kY=Z_l F$, where $Z_l=X_l\Gamma_k(FY_l)$. From Lemma~\ref{basklh12},
we have
$$
\|Z_l\|_{2, \mathcal{P}}\leq \widetilde{\alpha}_{k+1}\|X_l\|_{2, \mathcal{P}}\|Y_l\|_{2, \mathcal{P}}\leq
\widetilde{\alpha}_{k+1}\|X\|_{\mathcal{M}}\|Y\|_{\mathcal{M}}.
$$
Next, let $X=FX_r$ and $Y=FY_r$ with $X_r$ and $Y_r\in\mathfrak{S}_2(\mathcal{H}, \mathcal{P})$. Then $X\Gamma_kY=FZ_r$,
where $Z_r=X_r\Gamma_k(FY_r)$. Using Lemma~\ref{basklh12} once again, we obtain
$$
\|Z_r\|_{2, \mathcal{P}}\leq \widetilde{\alpha}_{k+1}\|X_r\|_{2, \mathcal{P}}\|Y_r\|_{2, \mathcal{P}}\leq
\widetilde{\alpha}_{k+1}\|X\|_{\mathcal{M}}\|Y\|_{\mathcal{M}}.
$$
The desired estimate for the norm of $(\Gamma_kX)Y$ is obtained in the same way, and the lemma is proved.
\epf

Theorem~\ref{baskth14} below is the main result of this section. It follows from Theorems~\ref{baskth6} and \ref{baskth12}, and Lemmata~\ref{basklh12} and \ref{basklh13}. The
property $\lim\limits_{k\to\infty}\widetilde{\alpha}_k=0$ ensures existence of $k\ge 0$ such that condition (\ref{bask12}) in
Theorem~\ref{baskth6}.

\bt\label{baskth14}
Let $m\in\ZZ_+$ be as in Theorem \ref{baskthsimilar} and
 $k\ge m$ be such that
$$
4\widetilde{\alpha}_{k+1}\|\widetilde{V}\|_{\mathcal{M}}<1.
$$
Then the operator $L_0-\widetilde{V}$ is similar to the operator $L_0-J_kX_*=L_0-V_0$, where  $X_*\in\mathcal{M}$
is the unique fixed point of the non-linear function $\Phi$ in (\ref{bask13}), where the transforms $J_k$ and $\Gamma_k$,
are defined by (\ref{Jm}) and (\ref{Gm}), and $B=\widetilde{V}$. Moreover, the operator $V_0$ is an orthogonal direct sum
$$
V_0=V_{0(k)}\oplus\left(\bigoplus_{|j|>k}V_{0,j}\right)
$$
with respect to the decomposition of the space $\mathcal{H}=L^2$ given by
$$
\mathcal{H}=\mathcal{H}_{(k)}\oplus\left(\bigoplus_{|j|>k}\mathcal{H}_j\right).
$$
In the above formula, $\mathcal{H}_{(k)}=\mathrm{Im}\,P_{(k)}$, $\mathcal{H}_j=\mathrm{Im}\,P_j$, $|j|>k$, where $P_{(k)}$ and
$P_j$, $|j|>k$, are the spectral projections of the operator $L_0$ defined by (\ref{bask10'}).
The similarity transform of  $L_0-\widetilde{V}$
into $L_0-V_0$ is the operator $I+\Gamma_kX_*$.
\et

\section{Proofs of the main results 
}\label{sect5}

Theorem~\ref{baskth3} follows from Theorems~\ref{baskthsimilar} and \ref{baskth14}. We remark that the similarity transform $U$ in Theorem~\ref{baskth3} has the form
\begeq\label{baskNew26}
U=U_{m,k}=(I+\Gamma_mV)(I+\Gamma_kX_*)=I+W_{m,k},
\eq
where $W_{m,k}=\Gamma_mV+\Gamma_kX_*+(\Gamma_mV)(\Gamma_kX_*)\in\mathfrak{S}_2(\mathcal{H}, \mathcal{P})$, and the numbers $k\ge m\ge 0$ can be determined from
 Theorems~\ref{baskthsimilar} and \ref{baskth14}.

 The result below follows from Theorem~\ref{baskth14} and Lemma~\ref{basklh1}.
\bt\label{baskth15}
The spectrum of the operator $L$ coincides with the spectrum of the operator
$$
L_0-V_0=L_0-P_{(k)}X_*P_{(k)}-\left(\bigoplus_{|j|>k}P_jX_*P_j\right).
$$
Moreover, we have
$$
\sigma(L)=\sigma(L_{(k)})\cup\sigma(L_j)=\sigma(L_{(k)})\cup\left(\bigcup_{|j|>k}\left(\frac{2\pi j}{\omega}+
\sigma(P_j{X_*}|{\mathcal{H}_j})\right)\right),
$$
where $L_{(k)}$ is the restriction of the operator $L_0-P_{(k)}X_*P_{(k)}$ to $\mathcal{H}_{(k)}=\mathrm{Im}\,P_{(k)}$ and the operators
$L_j$ are the restrictions of $L_0-P_jX_*P_j$ to $\mathrm{Im}\,P_j$, $|j|>k$.
\et

Theorem~\ref{baskth4} follows immediately from Theorem~\ref{baskth15}, since
$X_*\in\mathfrak{S}_2(\mathcal{H}, \mathcal{P})$ yields
 $\sum\limits_{|n|>k}\|P_nX_*P_n\|^2<\infty$.

Next, we estimate the spectral projections.

\bpf[Proof of Theorem~\ref{baskth5}]
Assume the conditions of Theorem~\ref{baskth3}.
Fix $k$, $m\in\mathbb{Z}_+$ that satisfy conditions of Theorems~\ref{baskthsimilar} and \ref{baskth14}. As usually, let
 ${P}_n=P(\{\lambda_n\}, L_0)$, $\lambda_n=\{\frac{i2\pi n}{\omega}\}$, $n\in\mathbb{Z}$, and
${P}_{(k)}=\sum\limits_{|n|\le k}{P}_n$ be the spectral projections for the operator $L_0$. Similarly, let 
$\widetilde{{P}}_{(k)}=U_{m,k}{P}_{(k)}U_{m,k}^{-1}$ and $\widetilde{{P}}_n=U_{m,k}{P}_nU_{m,k}^{-1}$ be the spectral projections for the operator $L$, as follows from
Lemma~\ref{basklh1}. Projections ${P}(\Omega)$ and $\widetilde{{P}}(\Omega)$ are defined by
$$
{P}(\Omega)=\sum_{n\in\Omega}{P}_n, \quad
\widetilde{{P}}(\Omega)=\sum_{n\in\Omega}\widetilde{{P}}_n=\sum_{n\in\Omega}U_{m,k}{P}_nU_{m,k}^{-1},
\quad \Omega=\mathbb{Z}\setminus\{-k, \dots, k\}.
$$
Obviously, $P(\Omega)$ is the spectral projection that corresponds to the spectral subset
$\bigcup\limits_{n\in\Omega}\{\lambda_n\}$ of the operator $L_0$.

Given $X\in\mathfrak{S}_2(\mathcal{H}, \mathcal{P})$, we let
$$
\alpha(\Omega, X)=\max\limits_{n\in\Omega}\alpha_n(X), \quad \Omega\subset\mathbb{Z},
$$
where the sequence $\alpha_n(X)$, $n\in\mathbb{Z}$, is given by (\ref{baskalpha}).

Let $X\in\mathcal{M}$, so that $X=X_l F_X=F_XX_r$, where $X_l$, $X_r\in\mathfrak{S}_2(\mathcal{H}, \mathcal{P})$. We have
\begin{equation*}
\bs
\|{P}(\Omega)X\|_{2, \mathcal{P}} & = \|{P}(\Omega)F_XX_r\|_{2, \mathcal{P}}
\\ & =
\Big\|\Big(\sum_{n\in\Omega}\alpha_n(X){P}_n\Big)X_r\Big\|_{2, \mathcal{P}}\le \alpha(\Omega, X)\|X\|_{2, \mathcal{P}}
\end{split}
\end{equation*}
and similarly for $X\mathcal{P}(\Omega)$. It follows that
$$
\max\{\|{P}(\Omega)X\|_{2, \mathcal{P}}, \|X{P}(\Omega)\|_{2, \mathcal{P}}\}\le  \|X\|_{\mathcal{M}}\alpha(\Omega, X).
$$
Next, we estimate $\|{P}(\Omega)\Gamma_kX\|_{2, \mathcal{P}}$ and $\|\Gamma_kX{P}(\Omega)\|_{2,
\mathcal{P}}$. We get
\begin{flalign*}
\|{P}(\Omega)\Gamma_kX\|_{2, \mathcal{P}}^2 &= \sum\limits_{\substack{\ell\in\Omega, j\in\mathbb{Z} \\ \ell\ne j}}
\frac{\|{P}_\ell{P}(\Omega)X{P}_j\|_{2, \mathcal{P}}^2}{|\lambda_\ell-\lambda_j|^2}\le 
\Big(\frac{\omega}{2\pi}\Big)^2\sum_{\ell\in\Omega}\|{P}_\ell X\|_{2, \mathcal{P}}^2 \\
&\le  \Big(\frac{\omega}{2\pi}\Big)^2\alpha^4(\Omega, X)\|X\|_{2, \mathcal{P}}^2.
\end{flalign*}
Therefore,
$$
\|{P}(\Omega)\Gamma_kX\|_{2, \mathcal{P}}\le  \frac{\omega}{2\pi}\alpha^2(\Omega, X)\|X\|_{2, \mathcal{P}},
$$
and similarly for $\|\Gamma_kX{P}(\Omega)\|_{2, \mathcal{P}}$. We note that the definition of the sequence $\alpha$ and the space $\mathcal{M}$ imply that
$\|F{P}(\Omega)\|=\|{P}(\Omega)F\|=\alpha(\Omega, B)$. Therfore, $\alpha(\Omega, X)\le \alpha(\Omega, B)\|X\|_{\mathcal{M}}$, $X\in\mathcal{M}$.

The two properties below also follow from the definition of $(\alpha_n(X))_n$.
\ben
\item If $X=\sum\limits_{\ell\ge  1}X_\ell$,  $X_\ell\in\mathfrak{S}_2(\mathcal{H}, \mathcal{P})$, then the series converges absolutely and
$$
\alpha_n(X)\|X\|_{2, \mathcal{P}}\le  \sum_{\ell\ge  1}\alpha_n(X_\ell)\|X_\ell\|_{2, \mathcal{P}}.
$$

\item If $X=X_1\cdot\ldots\cdot X_\ell$, $X_j\in\mathfrak{S}_2(\mathcal{H}, \mathcal{P})$, $1\le  j\le  \ell$,
then
$$
\alpha_n(X)\|X\|_{2, \mathcal{P}}\le  (\alpha_n(X_1)+\ldots+\alpha_n(X_\ell))\bigcap_{j=1}^\ell\|X_j\|_{2, \mathcal{P}}.
$$
\een

Next we estimate $\widetilde{{P}}(\Omega)-{P}(\Omega)$, where $\Omega\subset\mathbb{Z}
\setminus\{-k, \ldots, k\}$.
Recall that $U$ in Theorem~\ref{baskth3} is given by 
$U_{m,k}=I+W_{m,k}$, where $W_{m,k}=\Gamma_mV+\Gamma_kX_*+(\Gamma_mV)(\Gamma_kX_*)$. From Lemma~\ref{basklh1},
we get
$$
\widetilde{{P}}(\Omega)-{P}(\Omega)=(W_{m,k}{P}(\Omega)-{P}(\Omega)W_{m,k})(I+W_{m,k})^{-1}.
$$
Then
\begin{flalign*}
\|{P}(\Omega)W_{m,k}\|_{2, \mathcal{P}} &\le  \| {P}(\Omega)(\Gamma_mV)\|_{2, \mathcal{P}}+
\|{P}(\Omega)(\Gamma_kX_*)\|_{2, \mathcal{P}}  \\
&+ \|{P}(\Omega)(\Gamma_mV)(\Gamma_kX_*)\|_{2, \mathcal{P}} \\&\le 
\|{P}(\Omega)(\Gamma V)\|_{2, \mathcal{P}}+ \|{P}(\Omega)(\Gamma X_*)\|_{2, \mathcal{P}}  \\
&+ \|{P}(\Omega)(\Gamma V)(\Gamma X_*)\|_{2, \mathcal{P}}
\\& \le  C_1(\alpha(\Omega, \Gamma V) +
\alpha^2(\Omega, X_*))\le  C_2(\alpha(\Omega, \Gamma V) + \alpha^2(\Omega, V)),
\end{flalign*}
where the constants $C_1$ and $C_2$ are independent of $\Omega$. A similar estimate holds for
$\|W_{m,k}{P}(\Omega)\|_{2, \mathcal{P}}$.

The operator $(I+W_{m,k})^{-1}$ can be written as $(I+W_{m,k})^{-1}=I+\sum\limits_{j=1}^\infty(-1)^jW_{m,k}^j$, which yields
$$
\|(I+W_{m,k})^{-1}-I\|_{2, \mathcal{P}}\le  \frac{\|W_{m,k}\|_{2, \mathcal{P}}}{1-\|W_{m,k}\|_{2, \mathcal{P}}}.
$$
Collecting the above estimates together, we get
\begeq\label{bask26}
\|\widetilde{{P}}(\Omega)-{P}(\Omega)\|_{2, \mathcal{P}}\le  C_3(\alpha(\Omega, \Gamma V)+\alpha^2(\Omega, V)),
\eq
where $C_3$ is independent of $\Omega$.

Since the sequence $(\a_n(X))_n$ vanishes at infinity,
the assertion of Theorem~\ref{baskth5} immediately follows from (\ref{bask26}), if we let
 $\Omega=\{n\in\mathbb{Z}, |n|>N\}$ for sufficiently large $N\in\NN$.
\epf


To construct the group $T: \mathbb{R}\to B(L^2)$ generated by the operator $L: W_2^1\subset L^2\to L^2$, we use Theorem~\ref{baskth3} and its more detailed version~---
Theorem~\ref{baskth14}. Together with Theorem \ref{baskth7} and Lemma~\ref{baskuskova_lh6}, they imply that each operator $T(t)$, $t\in\mathbb{R}$, is a
$U$-orthogonal ($U=U_{m,k}$) direct sum of the form
\begin{equation}\label{bask27}
T(t)=U\left(\bigoplus_{j=-\infty}^{-k-1}e^{t(\frac{i2\pi j}{\omega}I_j+V_{0,j})}\oplus
e^{tV_{0(k)}}\oplus\left(\bigoplus_{j=k+1}^\infty e^{t(\frac{i2\pi j}{\omega}I_j+V_{0,j})}\right)\right)U^{-1}, 
\end{equation}
with respect to the $U$-orthogonal decomposition
$$
L^2=\mathcal{H}=\bigoplus_{j=-\infty}^{-k-1}U\mathcal{H}_j\oplus U\mathcal{H}_{(k)}\oplus
\left(\bigoplus_{j=k+1}^\infty U\mathcal{H}_j\right).
$$
The numbers $m$, $k\in\mathbb{Z}_+$, were defined in Theorems \ref{baskthsimilar} and \ref{baskth14}, respectively. The operators $V_{0(k)}$ and $V_{0,j}$ were defined in Theorem 
\ref{baskth14} as well.
The operator $U_{m,k}$ has the form (\ref{baskNew26}), so that $U_{m,k}=I+W_{m,k}$, where
$W_{m,k}\in\mathfrak{S}_2(\mathcal{H}, \mathcal{P})$.

Thus, the group $T: \mathbb{R}\to B(L^2)$ generated by $L=L_0-V$ has been constructed.
The assertions of Theorem~\ref{baskth1} regarding mild and classical solutions of Problem (\ref{bask1}) follow from the general theory of operator groups (see~\cite{DS58}, \cite{EN00}, \cite{HP57}). Formula (\ref{bask27}) yields Theorem~\ref{baskth6'} with $B_j=V_{0,j}$, $|j|\ge  k+1$, and
$B_{(k)}=V_{0(k)}$. From Theorems~\ref{baskth3} and \ref{baskth14} we get
$\sum\limits_{|j|\ge  k+1}\|V_{0,j}\|^2<\infty$.

Finally, Theorem~\ref{baskth8'} is obtained from (\ref{bask27}) and Parseval's identity.

\section{Examples}\label{sect7}

In this section, we present a concrete simple example that illustrates the method of similar operators and our main results.

We make the simplest possible choice of a non-trivial example. In particular, we let $H = \CC$, $\omega = 2\pi$, and choose a constant potential $\mathcal V (s) \equiv \frac1c$ with a sufficiently large $c\ge 1$. We then have $\mathcal H\simeq \ell^2(\mathbb Z)$, so that
\begeq\nonumber
L_0 \sim i\left(
\begin{array}{ccccccc}
\ddots  & \ddots & \ddots & \vdots & \iddots & \iddots & \iddots \\
\cdots & -2 & 0 & 0 &0  & 0  & \cdots \\
\cdots & 0  & -1 &0 & 0 & 0  & \cdots \\
\cdots & 0  &  0 &0 & 0 & 0  & \cdots \\
\cdots & 0  &  0 &0 &1  & 0   &\cdots \\
\cdots & 0  &  0 &0 &0  & 2   &\cdots \\
\iddots & \iddots & \iddots & \vdots & \ddots& \ddots & \ddots 
\end{array}
\right),
\quad (L_0)_{jk} = ij\delta_{j-k},
\eq 
where $\delta_j$ is the usual Kronecker delta, and
\begeq\nonumber
V \sim\frac1c\left(
\begin{array}{ccccccc}
\ddots  & \ddots & \ddots & \vdots & \iddots & \iddots & \iddots \\
\cdots & 0  & 0 & 0 &0  & 1  & \cdots \\
\cdots & 0  & 0 &0 & 1 & 0  & \cdots \\
\cdots & 0  & 0 &1 & 0 & 0  & \cdots \\
\cdots & 0  & 1 &0 &0  & 0   &\cdots \\
\cdots & 1  & 0 &0 &0  & 0   &\cdots \\
\iddots & \iddots & \iddots & \vdots & \ddots& \ddots & \ddots 
\end{array}
\right),
\quad (V)_{jk} = \frac1c\delta_{j+k}.
\eq 

We can compute the spectrum of $L_0 - V$ immediately by looking at $2\times 2$ matrices comprised of the rows and columns numbered $j$ and $-j$, $j\neq 0$:
\begeq\label{exspec}
\s(L) = \s(L_0-V) = \left\{-\frac1c\right\}\cup\left\{ij\sqrt{1-\frac1{c^2j^2}}: j\in\ZZ\setminus\{0\}\right\}.
\eq

To illustrate the set-up of our method we compute the following two matrices:
\begeq\nonumber
JV \sim\frac1c\left(
\begin{array}{ccccccc}
\ddots  & \ddots & \ddots & \vdots & \iddots & \iddots & \iddots \\
\cdots & 0  & 0 & 0 &0  & 0  & \cdots \\
\cdots & 0  & 0 &0 & 0 & 0  & \cdots \\
\cdots & 0  & 0 &1 & 0 & 0  & \cdots \\
\cdots & 0  & 0 &0 &0  & 0   &\cdots \\
\cdots & 0  & 0 &0 &0  & 0   &\cdots \\
\iddots & \iddots & \iddots & \vdots & \ddots& \ddots & \ddots 
\end{array}
\right),
\quad (JV)_{jk} = \frac1c\delta_{j}\delta_{k};
\eq

\begeq\nonumber
\Gamma V \sim\frac1{2ci}\left(
\begin{array}{ccccccc}
\ddots  & \ddots & \ddots & \vdots & \iddots & \iddots & \iddots \\
\cdots & 0  & 0 & 0 &0  & -\frac12  & \cdots \\
\cdots & 0  & 0 &0 & -1 & 0  & \cdots \\
\cdots & 0  & 0 &0 & 0 & 0  & \cdots \\
\cdots & 0  & 1 &0 &0  & 0   &\cdots \\
\cdots & \frac12  & 0 &0 &0  & 0   &\cdots \\
\iddots & \iddots & \iddots & \vdots & \ddots& \ddots & \ddots 
\end{array}
\right),
\quad (\Gamma V)_{jk} = \frac{\delta_{j+k}}{2jci}, j\neq 0.
\eq

A direct computation shows  $L_0\Gamma V - (\Gamma V) L_0 = V - JV$ so that \eqref{bask11'} is satisfied.

To illustrate the first similarity transform of our method we compute $U = (I+\Gamma V)^{-1}$ by looking at the same kind of $2\times 2$ submatrices as before
\eqref{exspec}. We get

\begeq\nonumber
U \sim\left(
\begin{array}{ccccccc}
\ddots  & \ddots & \ddots & \vdots & \iddots & \iddots & \iddots \\
\cdots & \frac{16c^2}{16c^2-1} & 0 & 0 &0  & - \frac{4ci}{16c^2-1}  & \cdots \\
\cdots & 0  & \frac{4c^2}{4c^2-1} &0 &  - \frac{2ci}{4c^2-1} & 0  & \cdots \\
\cdots & 0  &  0 &1 & 0 & 0  & \cdots \\
\cdots & 0  &   \frac{2ci}{4c^2-1} &0 & \frac{4c^2}{4c^2-1}  & 0   &\cdots \\
\cdots &  \frac{4ci}{16c^2-1}  &  0 &0 &0  &  \frac{16c^2}{16c^2-1}   &\cdots \\
\iddots & \iddots & \iddots & \vdots & \ddots& \ddots & \ddots 
\end{array}
\right),
\eq
$$
\quad (U)_{jk} = \frac{4j^2c^2\delta_{j-k}}{4j^2c^2-1}+ \frac{2jci\delta_{j+k}}{4j^2c^2-1},\ j\neq 0, \ (U)_{00} = 1.
$$
Recall that we chose $c$ sufficiently large so that the above matrix can also be computed via a geometric series:
$$U = (I+\Gamma V)^{-1} = \sum_{n = 0}^\infty (-\Gamma V)^n = \sum_{n = 0}^\infty (\Gamma V)^{2n} -
\Gamma V\sum_{n = 0}^\infty (\Gamma V)^{2n}.$$ 
The choice of $c$ also guarantees that we do not need to deal with coarser resolutions of the identity neither in the first nor in the second similarity transform
(we have $m = k =0$ in Theorem \ref{baskth14}).
Since
\begeq\nonumber
V\Gamma V \sim \frac1{2ic^2}\left(
\begin{array}{ccccccc}
\ddots  & \ddots & \ddots & \vdots & \iddots & \iddots & \iddots \\
\cdots & \frac12 & 0 & 0 &0  & 0  & \cdots \\
\cdots & 0  & 1 &0 & 0 & 0  & \cdots \\
\cdots & 0  &  0 &0 & 0 & 0  & \cdots \\
\cdots & 0  &  0 &0 &-1  & 0   &\cdots \\
\cdots & 0  &  0 &0 &0  & -\frac12   &\cdots \\
\iddots & \iddots & \iddots & \vdots & \ddots& \ddots & \ddots 
\end{array}
\right),
\quad (V\Gamma V)_{jk} = \frac{i\delta_{j-k}}{2jc^2}, j\neq 0,
\eq 
and, from \eqref{basktildev},
$
\widetilde V = JV + U(V\Gamma V - (\Gamma V)JV) = JV + UV\Gamma V,
$
we get
\begeq\nonumber
\widetilde V \sim 
\left(
\begin{array}{ccccccc}
\ddots  & \ddots & \ddots & \vdots & \iddots & \iddots & \iddots \\
\cdots &  -\frac{4i}{16c^2-1}  & 0 & 0 &0  &  \frac{1/c}{16c^2-1}  & \cdots \\
\cdots & 0  & -\frac{2i}{4c^2-1} &0 & \frac{1/c}{4c^2-1} & 0  & \cdots \\
\cdots & 0  & 0 &\frac1c & 0 & 0  & \cdots \\
\cdots & 0  & \frac{1/c}{4c^2-1} &0 &\frac{2i}{4c^2-1}  & 0   &\cdots \\
\cdots &  \frac{1/c}{16c^2-1}  & 0 &0 &0  & \frac{4i}{16c^2-1}   &\cdots \\
\iddots & \iddots & \iddots & \vdots & \ddots& \ddots & \ddots 
\end{array}
\right),
\eq
\[
\quad (\widetilde V)_{jk} = \frac{\delta_{j+k}}{c(4j^2c^2-1)}+  \frac{2ji\delta_{j-k}}{4j^2c^2-1},\ j\neq 0, (\widetilde V)_{00} =\frac1{c}.
\]
Thus, $\widetilde {V} \in {\mathfrak S}_2(\HH)$, which was the goal of the first similarity transform. Observe also that the $2\times2$ submatrices of $L_0-\widetilde V$ are
\begeq\nonumber
\frac1{4c^2j^2-1}\left(
\begin{array}{cc}
-ij(4c^2j^2-3)  & 1/c  \\
1/c & ij(4c^2j^2-3) 
\end{array}
\right),
\quad  j\neq 0.
\eq
After a short direct computation, we see that $L_0-V$ and $L_0-\widetilde V$ do indeed have the same spectrum \eqref{exspec}.

Next, we would like to illustrate our second similarity transform by finding the fixed point $X_*$ of the function \eqref{bask13} with $B = \widetilde V$.
To simplify the notation, we let $v_j = \frac 1{c(4j^2c^2-1)}$ so that
\[
\quad (\widetilde V)_{jk} = v_j\delta_{j+k} +  2cjiv_j\delta_{j-k},\ j\neq 0. 
\]
We 
also let $x_j = (X_*)_{jj}$ and $y_j = (X_*)_{j,-j}$. Observe that all other entries of the
matrix of $X_*$ are $0$, so that it suffices to look at the same kind of the $2\times 2$ submatrices as in the computation of the spectrum of $L_0-V$.
Thus, from $\Phi(X_*) = X_* = \widetilde V\Gamma X_* - (\Gamma X_*)J\widetilde V -(\Gamma X_*)J(\widetilde V\Gamma X_*) +\widetilde V$ 
we must have
\begeq\nonumber
\bs
\left(
\begin{array}{cc}
x_{-j}  & y_{-j}  \\
y_j & x_j 
\end{array}
\right) & =
v_j\left(
\begin{array}{cc}
-2cij  & 1  \\
1 & 2cij 
\end{array}
\right) \cdot
\frac1{2ij}
\left(
\begin{array}{cc}
0  & -y_{-j}  \\
y_j & 0 
\end{array}
\right)
\\ &-
\frac{1}{2ij}\left(
\begin{array}{cc}
0  & -y_{-j}  \\
y_j & 0 
\end{array}
\right) \cdot
v_j\left(
\begin{array}{cc}
-2cij  & 0  \\
0 & 2cij 
\end{array}
\right)
\\
&-
\frac{1}{2ij}\left(
\begin{array}{cc}
0  & -y_{-j}  \\
y_j & 0 
\end{array}
\right) \cdot 
v_j\left(
\begin{array}{cc}
0  & 1  \\
1 & 0 
\end{array}
\right) \cdot
\frac1{2ij}
\left(
\begin{array}{cc}
0  & -y_{-j}  \\
y_j & 0 
\end{array}
\right)
\\
&+v_j\left(
\begin{array}{cc}
-2cij  & 1  \\
1 & 2cij 
\end{array}
\right) 
\\
&=
\frac{v_j}{2ij}
\left[
\left(
\begin{array}{cc}
y_j& 2cijy_{-j}    \\
2cijy_{j}  & y_{-j} 
\end{array}
\right) -
\left(
\begin{array}{cc}
0& -2cijy_{-j}    \\
-2cijy_{j}  & 0 
\end{array}
\right)
\right]
\\
&+
\frac{v_j}{4j^2}
\left(
\begin{array}{cc}
0  & y_{-j}^2  \\
y_j^2 & 0 
\end{array}
\right)+v_j\left(
\begin{array}{cc}
-2cij  & 1  \\
1 & 2cij 
\end{array}
\right) 
\\
& = 
\frac{v_j}{2ij}
\left(
\begin{array}{cc}
y_j+4cj^2  & 4cijy_{-j} -\frac{y_{-j}^2}{2ij}+2ij \\
4cijy_{j}-\frac{y_{j}^2}{2ij} +2ij& y_{-j}-4cj^2 
\end{array}
\right) .
\end{split}
\eq
With a little bit of work, it follows that $y_{-j}=y_j$ and
\[
\left\{
\bs
x_j &= \frac{iv_j}{2j}y_{j} + 2cjiv_j\\
y_j &= 2cv_jy_j + \frac{v_j}{4j^2}y_j^2+v_j
\end{split}
\right.,\ j\neq 0.
\]
Therefore, we get 
\[
\left\{
\bs
x_j &= ij\left(1 - \sqrt{\left(1-2cv_j\right)^2-\frac{v_j^2}{j^2}}\right)\\
y_j &= -4cj^2+\frac{2j^2}{v_j}\left(1 - \sqrt{\left(1-2cv_j\right)^2-\frac{v_j^2}{j^2}}\right)
\end{split}
\right.,
\]
which, after simplification, implies
\[
x_j = ij \left(1- {\sqrt{1-\frac1{c^2j^2}}}\right) = \frac i{c^2j\left(1+ {\sqrt{1-\frac1{c^2j^2}}}\right)}, \ j\neq0.
\]
The above equality is once again consistent with \eqref{exspec} as 
it is immediate that   $x_0 = -\frac1c$. Recall that Theorem \ref{baskth6} yields that $L_0-\widetilde V$ is similar to $L_0 - JX_*$ and we have just confirmed
that, in this example, the two operators do, indeed, have the same spectrum. 
Observe also that $\{x_j\}$ is  a square summable sequence as was predicted by  Theorem \ref{baskth4}.
This highlights a major difference that is caused by the presence of the involution. In an analogous example without the involution,
we would be looking at the operator $L_0 - \frac1c I$ which has a very different spectral structure in terms of the perturbation.

Of course, in a more sophisticated example, one may be unable to find the fixed point $X_*$ explicitly. This underscores the importance of the estimates in Theorem \ref{baskth8'} and the equiconvergence result of Theorem \ref{baskth5}.

\section{Acknowledgment}

We thank the anonymous referee for their helpful suggestions to improve the manuscript.

\bibliographystyle{siam}
\bibliography{../refs}

\begin{thebibliography}{10}

\bibitem{ABHN11}
{\sc W.~Arendt, C.~J.~K. Batty, M.~Hieber, and F.~Neubrander}, {\em
  Vector-valued {L}aplace transforms and {C}auchy problems}, vol.~96 of
  Monographs in Mathematics, Birkh\"auser/Springer Basel AG, Basel, second~ed.,
  2011.

\bibitem{BD15}
{\sc A.~Baskakov and V.~Didenko}, {\em Spectral analysis of differential
  operators with unbounded periodic coefficients}, Differential Equations, 51
  (2015), pp.~325--341.

\bibitem{B78}
{\sc A.~G. Baskakov}, {\em Spectral tests for the almost periodicity of the
  solutions of functional equations}, Mat. Zametki, 24 (1978), pp.~195--206,
  301.
\newblock English translation: Math. Notes 24 (1978), no. 1--2, pp. 606--612
  (1979).

\bibitem{B83}
\leavevmode\vrule height 2pt depth -1.6pt width 23pt, {\em Methods of abstract
  harmonic analysis in the theory of perturbations of linear operators},
  Sibirsk. Mat. Zh., 24 (1983), pp.~21--39, 191.
\newblock English translation: Siberian Math. J. 24 (1983), no. 1, pp. 17--32.

\bibitem{B84u}
\leavevmode\vrule height 2pt depth -1.6pt width 23pt, {\em The
  {K}rylov-{B}ogolyubov substitution in the theory of perturbations of linear
  operators}, Ukrain. Mat. Zh., 36 (1984), pp.~606--611.
\newblock English translation: Ukrainian Math. J. 36 (1984), no. 5, 451--455.

\bibitem{B86}
\leavevmode\vrule height 2pt depth -1.6pt width 23pt, {\em A theorem on
  splitting of an operator and some related problems in the analytic theory of
  perturbations}, Izv. Akad. Nauk SSSR Ser. Mat., 50 (1986), pp.~435--457, 638.
\newblock English translation: Math. USSR-Izv. 28 (1987), no. 3, pp. 421--444.

\bibitem{B94}
\leavevmode\vrule height 2pt depth -1.6pt width 23pt, {\em Spectral analysis of
  perturbed non-quasi-analytic and spectral operators}, Izv. Ross. Akad. Nauk
  Ser. Mat., 58 (1994), pp.~3--32.
\newblock English translation: Russian Acad. Sci. Izv. Math. 45 (1995), no. 1,
  pp. 1--31.

\bibitem{B15}
\leavevmode\vrule height 2pt depth -1.6pt width 23pt, {\em Estimates for the
  {G}reen's function and parameters of exponential dichotomy of a hyperbolic
  operator semigroup and linear relations}, Mat. Sb., 206 (2015), pp.~23--62.
\newblock English translation: Sb. Math., 206 (2015), no. 8, pp. 1049--1086.

\bibitem{BDS11}
{\sc A.~G. Baskakov, A.~V. Derbushev, and A.~O. Shcherbakov}, {\em The method
  of similar operators in the spectral analysis of the nonselfadjoint {D}irac
  operator with nonsmooth potential}, Izv. Ross. Akad. Nauk Ser. Mat., 75
  (2011), pp.~3--28.
\newblock English translation: Izv. Math. 75 (2011), no. 3, 445--469.

\bibitem{BKKS17}
{\sc A.~G. Baskakov, L.~Y. Kabantsova, I.~D. Kostrub, and T.~I. Smagina}, {\em
  Linear differential operators and operator matrices of the second order},
  Differ. Equ., 53 (2017), pp.~8--17.
\newblock Translation of Differ. Uravn. {{\bf{5}}3} (2017), no. 1, 10--19.

\bibitem{BK13}
{\sc A.~G. Baskakov and I.~A. Krishtal}, {\em On completeness of spectral
  subspaces of linear relations and ordered pairs of linear operators}, J.
  Math. Anal. Appl., 407 (2013), pp.~157--178.

\bibitem{BK17}
\leavevmode\vrule height 2pt depth -1.6pt width 23pt, {\em Spectral analysis of
  abstract parabolic operators in homogeneous function spaces, ii},
  Mediterranean Journal of Mathematics, 14 (2017), p.~181.

\bibitem{BKR17}
{\sc A.~G. Baskakov, I.~A. Krishtal, and E.~Y. Romanova}, {\em Spectral
  analysis of a differential operator with an involution}, Journal of Evolution
  Equations, 17 (2017), pp.~669--684.

\bibitem{BP17}
{\sc A.~G. Baskakov and D.~M. Polyakov}, {\em The method of similar operators
  in spectral analysis for the {H}ill operator with nonsmooth potential}, Mat.
  Sb., 208 (2017), pp.~3--47.
\newblock English Translation: Sb. Math., 2017, 208(1): 1--43.

\bibitem{BKh11z}
{\sc M.~S. Burlutskaya and A.~P. Khromov}, {\em The {F}ourier method in a mixed
  problem for a first-order partial differential equation with involution}, Zh.
  Vychisl. Mat. Mat. Fiz., 51 (2011), pp.~2233--2246.
\newblock English translation: Comput. Math. Math. Phys. 51 (2011), no. 12,
  2102--2114.

\bibitem{BKh11}
\leavevmode\vrule height 2pt depth -1.6pt width 23pt, {\em Mixed problems for
  first-order hyperbolic equations with involution}, Dokl. Akad. Nauk, 441
  (2011), pp.~156--159.
\newblock English translation: Dokl. Math. 84 (2011), no. 3, 783--786.

\bibitem{BKh14}
\leavevmode\vrule height 2pt depth -1.6pt width 23pt, {\em
  Functional-differential operators with involution and {D}irac operators with
  periodic boundary conditions}, Dokl. Akad. Nauk, 454 (2014), pp.~15--17.
\newblock English translation: Dokl. Math. 89 (2014), no. 1, 8--10.

\bibitem{BKK12}
{\sc M.~S. Burlutskaya, V.~P. Kurdyumov, and A.~P. Khromov}, {\em Refined
  asymptotic formulas for the eigenvalues and eigenfunctions of the {D}irac
  system}, Dokl. Akad. Nauk, 443 (2012), pp.~414--417.
\newblock English translation: Dokl. Math. 85 (2012), no. 2, 240--242.

\bibitem{CL99}
{\sc C.~Chicone and Y.~Latushkin}, {\em Evolution semigroups in dynamical
  systems and differential equations}, vol.~70 of Mathematical Surveys and
  Monographs, American Mathematical Society, Providence, RI, 1999.

\bibitem{DS58}
{\sc N.~Dunford and J.~T. Schwartz}, {\em Linear {O}perators. {I}. {G}eneral
  {T}heory}, With the assistance of W. G. Bade and R. G. Bartle. Pure and
  Applied Mathematics, Vol. 7, Interscience Publishers, Inc., New York, 1958.

\bibitem{EN00}
{\sc K.-J. Engel and R.~Nagel}, {\em One-parameter semigroups for linear
  evolution equations}, vol.~194 of Graduate Texts in Mathematics,
  Springer-Verlag, New York, 2000.
\newblock With contributions by S. Brendle, M. Campiti, T. Hahn, G. Metafune,
  G. Nickel, D. Pallara, C. Perazzoli, A. Rhandi, S. Romanelli and R.
  Schnaubelt.

\bibitem{F65}
{\sc K.~O. Friedrichs}, {\em Lectures on advanced ordinary differential
  equations}, Notes by P. Berg, W. Hirsch, P. Treuenfels, Gordon and Breach
  Science Publishers, New York, 1965.

\bibitem{GKK96}
{\sc I.~Gohberg, M.~A. Kaashoek, and J.~Kos}, {\em Classification of linear
  time-varying difference equations under kinematic similarity}, Integral
  Equations Operator Theory, 25 (1996), pp.~445--480.

\bibitem{HP57}
{\sc E.~Hille and R.~S. Phillips}, {\em Functional analysis and semi-groups},
  American Mathematical Society Colloquium Publications, vol. 31, American
  Mathematical Society, Providence, R. I., 1957.
\newblock rev. ed.

\bibitem{KB61}
{\sc R.~E. Kalman and R.~S. Bucy}, {\em New results in linear filtering and
  prediction theory}, Trans. ASME Ser. D. J. Basic Engrg., 83 (1961),
  pp.~95--108.

\bibitem{K81}
{\sc A.~P. Khromov}, {\em Equiconvergence theorems for integro-differential and
  integral operators}, Mat. Sb. (N.S.), 114(156) (1981), pp.~378--405, 479.

\bibitem{LMS95}
{\sc Y.~Latushkin and S.~Montgomery-Smith}, {\em Evolutionary semigroups and
  {L}yapunov theorems in {B}anach spaces}, J. Funct. Anal., 127 (1995),
  pp.~173--197.

\bibitem{Lj56}
{\sc A.~M. Ljapunov}, {\em { {S}obranie sochineni\u\i . {T}om}\ {II}}, Izdat.
  Akad. Nauk SSSR, Moscow, 1956.

\bibitem{Mit04}
{\sc B.~Mityagin}, {\em Spectral expansions of one-dimensional periodic {D}irac
  operators}, Dyn. Partial Differ. Equ., 1 (2004), pp.~125--191.

\bibitem{N15}
{\sc M.~Niezabitowski}, {\em Kinematic similarity of the discrete linear
  time-varying systems}, in 2015 20th International Conference on Control
  Systems and Computer Science, May 2015, pp.~10--17.

\bibitem{P64}
{\sc V.~A. Pliss}, {\em Nelokalnye problemy teorii kolebanii}, Izdat.
  ``Nauka'', Moscow, 1964.

\bibitem{P65}
\leavevmode\vrule height 2pt depth -1.6pt width 23pt, {\em Families of periodic
  solutions of systems of differential equations of second order without
  dissipation}, Differencial'nye Uravnenija, 1 (1965), pp.~1428--1448.

\bibitem{P16}
{\sc D.~M. Polyakov}, {\em Spectral properties of an even-order differential
  operator}, Differ. Equ., 52 (2016), pp.~1098--1103.
\newblock Translation of Differ. Uravn. {{\bf{5}}2} (2016), no. 8, 1133--1137.

\bibitem{SS12}
{\sc M.~A. Sadybekov and A.~M. Sarsenbi}, {\em Criterion for the basis property
  of the eigenfunction system of a multiple differentiation operator with an
  involution}, Differ. Equ., 48 (2012), pp.~1112--1118.
\newblock Translation of Differ. Uravn. {{\bf{4}}8} (2012), no. 8, 1126--1132.

\bibitem{T01}
{\sc Y.~Tomilov}, {\em A resolvent approach to stability of operator
  semigroups}, J. Operator Theory, 46 (2001), pp.~63--98.

\bibitem{T65}
{\sc R.~E.~L. Turner}, {\em Perturbation of compact spectral operators}, Comm.
  Pure Appl. Math., 18 (1965), pp.~519--541.

\bibitem{U04}
{\sc N.~B. Uskova}, {\em On a result of {R}. {T}urner}, Mat. Zametki, 76
  (2004), pp.~905--917.
\newblock English translation: Math. Notes 76 (2004), no. 5-6, 844--854.

\bibitem{U15}
\leavevmode\vrule height 2pt depth -1.6pt width 23pt, {\em On spectral
  properties of {S}turm-{L}iouville operator with matrix potential}, Ufa Math.
  J., 7 (2015), pp.~84--94.

\bibitem{U16}
\leavevmode\vrule height 2pt depth -1.6pt width 23pt, {\em On the spectral
  properties of a second-order differential operator with a matrix potential},
  Differ. Equ., 52 (2016), pp.~557--567.
\newblock Translation of Differ. Uravn. {{\bf{5}}2} (2016), no. 5, 579--588.

\bibitem{V84}
{\sc V.~S. Vladimirov}, {\em Equations of mathematical physics}, ``Mir'',
  Moscow, 1984.
\newblock Translated from the Russian by Eugene Yankovsky [E. Yankovski{\u\i}].

\end{thebibliography}

\end {document}